\def\DATE{November 19, 2013}
\def\DATE{\relax}
\magnification=1100
\baselineskip=12.72pt
\voffset=.75in
\hoffset=.9in
\hsize=4.1in
\newdimen\hsizeGlobal
\hsizeGlobal=4.1in%
\vsize=7.05in
\parindent=.166666in
\pretolerance=500 \tolerance=1000 \brokenpenalty=5000

\footline={\vbox{\hsize=\hsizeGlobal\hfill{\rm\the\pageno}\hfill\llap{\sevenrm\DATE}}\hss}

\long\def\note#1{%
  \hfuzz=50pt%
  \vadjust{%
    \setbox1=\vtop{%
      \hsize 3cm\parindent=0pt\eightpoints\baselineskip=9pt%
      \rightskip=4mm plus 4mm\raggedright#1\hss%
      }%
    \hbox{\kern-4cm\smash{\box1}\hss\par}%
    }%
  \hfuzz=0pt
  }
\def\note#1{\relax}

\def\anote#1#2#3{\smash{\kern#1in{\raise#2in\hbox{#3}}}%
  \nointerlineskip}     
\def\anotePT#1#2#3{\smash{\kern#1pt{\raise#2pt\hbox{#3}}}%
  \nointerlineskip}     

\newcount\equanumber
\equanumber=0
\newcount\sectionnumber
\sectionnumber=0
\newcount\subsectionnumber
\subsectionnumber=0
\newcount\snumber  
\snumber=0

\def\section#1{%
  \subsectionnumber=0%
  \snumber=0%
  \equanumber=0%
  \advance\sectionnumber by 1%
  \noindent{\bf \the\sectionnumber .~#1.~}%
}%
\def\subsection#1{%
  \advance\subsectionnumber by 1%
  \snumber=0%
  \equanumber=0%
  \noindent{\bf \the\sectionnumber .\the\subsectionnumber .~#1.~}%
}%
\def\prevs{\the\sectionnumber .\the\subsectionnumber .\the\snumber }
\long\def\Definition#1{%
  \global\advance\snumber by 1%
  \bigskip
  \noindent{\bf Definition~\prevs .}%
  \quad{\it#1}%
}

\long\def\Corollary#1{%
  \global\advance\snumber by 1%
  \bigskip
  \noindent{\bf Corollary~\prevs .}%
  \quad{\it#1}%
}%
\long\def\Lemma#1{%
  \global\advance\snumber by 1%
  \bigskip
  \noindent{\bf Lemma~\prevs .}%
  \quad{\it#1}%
}%
\long\def\Notation#1{%
  \global\advance\snumber by 1%
  \bigskip
  \noindent{\bf Notation~\prevs .}%
  \quad{\it#1}%
}%
\def\Proof{\noindent{\bf Proof.~}}
\long\def\Proposition#1{%
  \advance\snumber by 1%
  \bigskip
  \noindent{\bf Proposition~\prevs .}%
  \quad{\it#1}%
}%
\long\def\Remark#1{%
  \bigskip
  \noindent{\bf Remark.~}#1%
}%
\long\def\remark#1{%
  \advance\snumber by1%
  \bigskip
  \noindent{\bf Remark~\prevs .}\quad#1%
}%
\long\def\Theorem#1{%
  \advance\snumber by 1%
  \bigskip
  \noindent{\bf Theorem~\prevs .}%
  \quad{\it#1}%
}%
\long\def\Statement#1{%
  \advance\snumber by 1%
  \bigskip
  \noindent{\bf Statement~\prevs .}%
  \quad{\it#1}%
}%
\def\ifundefined#1{\expandafter\ifx\csname#1\endcsname\relax}
\def\labeldef#1{\global\expandafter\edef\csname#1\endcsname{\prevs}}
\def\labelref#1{\expandafter\csname#1\endcsname}
\def\label#1{\ifundefined{#1}\labeldef{#1}\note{$<$#1$>$}\else\labelref{#1}\fi}

\def\preveq{(\the\sectionnumber .\the\subsectionnumber .\the\equanumber)}
\def\neq{\global\advance\equanumber by 1\eqno{\preveq}}

\def\ifundefined#1{\expandafter\ifx\csname#1\endcsname\relax}

\def\equadef#1{\global\advance\equanumber by 1%
  \global\expandafter\edef\csname#1\endcsname{\preveq}%
  \setbox1=\hbox{\rm\hskip .1in[#1]}\dp1=0pt\ht1=0pt\wd1=0pt%
  \preveq\box1}
\def\equadef#1{\global\advance\equanumber by 1%
  \global\expandafter\edef\csname#1\endcsname{\preveq}%
  \preveq}

\def\equaref#1{\expandafter\csname#1\endcsname}

\def\equa#1{%
  \ifundefined{#1}%
    \equadef{#1}%
  \else\equaref{#1}\fi}

\font\eightrm=cmr8%
\font\sixrm=cmr6%

\font\eightsl=cmsl8%

\font\eightbf=cmb8%

\font\eighti=cmmi8%
\font\sixi=cmmi6%

\font\eightsy=cmsy8%
\font\sixsy=cmsy6%

\font\eightex=cmex8%
\font\sixex=cmex6%
\font\fiveex=cmex5%

\font\eightit=cmti8%

\font\eighttt=cmtt8%

\font\tenbb=msbm10%
\font\eightbb=msbm8%
\font\sevenbb=msbm7%
\font\sixbb=msbm6%
\font\fivebb=msbm5%
\newfam\bbfam  \textfont\bbfam=\tenbb  \scriptfont\bbfam=\sevenbb  \scriptscriptfont\bbfam=\fivebb%

\font\tenbbm=bbm10

\font\tencmssi=cmssi10%
\font\sevencmssi=cmssi7%
\font\fivecmssi=cmssi5%
\newfam\ssfam  \textfont\ssfam=\tencmssi  \scriptfont\ssfam=\sevencmssi  \scriptscriptfont\ssfam=\fivecmssi%

\font\tenfrak=cmfrak10%
\font\eightfrak=cmfrak8%
\font\sevenfrak=cmfrak7%
\font\sixfrak=cmfrak6%
\font\fivefrak=cmfrak5%
\newfam\frakfam  \textfont\frakfam=\tenfrak  \scriptfont\frakfam=\sevenfrak  \scriptscriptfont\frakfam=\fivefrak%
\def\frak{\fam\frakfam\tenfrak}%

\font\tenmsam=msam10%
\font\eightmsam=msam8%
\font\sevenmsam=msam7%
\font\sixmsam=msam6%
\font\fivemsam=msam5%

\def\bb{\fam\bbfam\tenbb}%

\def\hexdigit#1{\ifnum#1<10 \number#1\else%
  \ifnum#1=10 A\else\ifnum#1=11 B\else\ifnum#1=12 C\else%
  \ifnum#1=13 D\else\ifnum#1=14 E\else\ifnum#1=15 F\fi%
  \fi\fi\fi\fi\fi\fi}
\newfam\msamfam  \textfont\msamfam=\tenmsam  \scriptfont\msamfam=\sevenmsam  \scriptscriptfont\msamfam=\fivemsam%
\def\msam{\msamfam\tenmsam}%
\mathchardef\leq"3\hexdigit\msamfam 36%
\mathchardef\geq"3\hexdigit\msamfam 3E%

\font\tentt=cmtt11%
\font\seventt=cmtt9%
\textfont\ttfam=\tentt
\scriptfont7=\seventt%
\def\tt{\fam\ttfam\tentt}%

\def\eightpoints{%
\def\rm{\fam0\eightrm}%
\textfont0=\eightrm   \scriptfont0=\sixrm   \scriptscriptfont0=\fiverm%
\textfont1=\eighti    \scriptfont1=\sixi    \scriptscriptfont1=\fivei%
\textfont2=\eightsy   \scriptfont2=\sixsy   \scriptscriptfont2=\fivesy%
\textfont3=\eightex   \scriptfont3=\sixex   \scriptscriptfont3=\fiveex%
\textfont\itfam=\eightit  \def\it{\fam\itfam\eightit}%
\textfont\slfam=\eightsl  \def\sl{\fam\slfam\eightsl}%
\textfont\ttfam=\eighttt  \def\tt{\fam\ttfam\eighttt}%
\textfont\bffam=\eightbf  \def\bf{\fam\bffam\eightbf}%

\textfont\frakfam=\eightfrak  \scriptfont\frakfam=\sixfrak \scriptscriptfont\frakfam=\fivefrak  \def\frak{\fam\frakfam\eightfrak}%
\textfont\bbfam=\eightbb      \scriptfont\bbfam=\sixbb     \scriptscriptfont\bbfam=\fivebb      \def\bb{\fam\bbfam\eightbb}%
\textfont\msamfam=\eightmsam  \scriptfont\msamfam=\sixmsam \scriptscriptfont\msamfam=\fivemsam  \def\msam{\msamfam\eightmsam}

\rm%
}

\def\poorBold#1{\setbox1=\hbox{#1}\wd1=0pt\copy1\hskip0.25pt\box1\hskip0.25pt#1}%

\mathchardef\lsim"3\hexdigit\msamfam 2E%
\mathchardef\gsim"3\hexdigit\msamfam 26%

\long\def\algorithm#1{{\parindent=0pt\tt\par#1\par}}

\def\ds{\displaystyle}
\long\def\DoNotPrint#1{\relax}
\def\Deg{{\rm Deg}}

\def\fixedref#1{#1\note{fixedref$\{$#1$\}$}}

\def\limn{\lim_{n\to\infty}}

\def\mod{\mathop{\rm mod}}
\def\ord{\mathop{\rm ord}}

\def\qed{{\vrule height .9ex width .8ex depth -.1ex}}
\def\ss{\scriptstyle}

\def\boc{\note{{\bf BoC}\hskip-11pt\setbox1=\hbox{$\Bigg\downarrow$}%
         \dp1=0pt\ht1=0pt\ht1=0pt\leavevmode\raise -20pt\box1}}
\def\eoc{\note{{\bf EoC}\hskip-11pt\setbox1=\hbox{$\Bigg\uparrow$}%
         \dp1=0pt\ht1=0pt\ht1=0pt\leavevmode\raise 20pt\box1}}

\def\One{\hbox{\tenbbm 1}}

\def\og{\overline{g}}

\def\ug{\underline{g}}

\def\calA{{\cal A}}
\def\calB{{\cal B}}
\def\calC{{\cal C}}

\def\calO{{\cal O}}

\def\calQ{{\cal Q}}

\def\CC{{\bb C}}

\def\MM{{\bb M}\kern .4pt}
\def\NN{{\bb N}\kern .5pt}

\def\ZZ{{\bb Z}}

%
%
%
%
\def\uncatcodespecials 
    {\def\do##1{\catcode`##1=12}\dospecials}%
{\catcode`\^^I=\active \gdef^^I{\ \ \ \ }
 \catcode`\`=\active\gdef`{\relax\lq}}
\def\setupverbatim 
    {\parindent=0pt\tt %
     \spaceskip=0pt \xspaceskip=0pt 
     \catcode`\^^I=\active %
     \catcode`\`=\active %
     \def\par{\leavevmode\endgraf}
     \obeylines \uncatcodespecials \obeyspaces %
     }%
{\obeyspaces \global\let =\ }
%
%
%
%
%
%

\pageno=1

\centerline{\bf ON \poorBold{$q$}-ALGEBRAIC EQUATIONS}
\centerline{\bf AND THEIR POWER SERIES SOLUTIONS}

\bigskip
 
\centerline{Ph.\ Barbe$^{(1)}$ and W.P.\ McCormick$^{(2)}$}
\centerline{${}^{(1)}$CNRS {\sevenrm(UMR {\eightrm 8088})}, ${}^{(2)}$University of Georgia}

 
{\narrower
\baselineskip=9pt\parindent=0pt\eightpoints

\bigskip

{\bf Abstract.} We study the existence of formal power series solutions 
to $q$-algebraic equations. When a solution exists, we give a sufficient
condition on the equation for this solution to have a positive
radius of convergence. We emphasize on the case where the solution is 
divergent, giving a sharp
estimate on the growth of the coefficients. As a consequence,
we obtain a bound on the $q$-Gevrey order of the formal solution, which is
optimal in some cases. Various examples illustrate our main results.

\bigskip


\bigskip
 

}

\bigskip

\def\prevs{\the\sectionnumber.\the\snumber }
\def\preveq{(\the\sectionnumber.\the\equanumber)}

\section{Introduction}
The purpose of this paper is to develop a theory of $q$-algebraic equations
and their power series solutions, meaning the following: a formal power series 
$f(z)=\sum_{n\geq 0} f_n z^n$ with complex coefficients $(f_n)$ is $q$-algebraic
if there exists a complex polynomial of several variables 
$P(z,z_0,z_1,\ldots,z_k)$ and a complex number $q$ such that
$$
  P\bigl(z,f(z),f(qz),\ldots,f(q^kz)\bigr)=0 \, .
  \eqno{\equa{qAlgebraic}}
$$
Since the coefficients of the polynomial are free parameters in this paper, 
they may in fact depend on $q$ in any fashion that one wishes. These
$q$-algebraic equations occur in various areas, including combinatorics,
dynamical systems, knot theory, and mathematical physics.

There is no loss of generality in assuming that the polynomial $P$ is 
irreducible. However, in some applications, such as in combinatorics, one may
derive a $q$-algebraic equation for a generating function, and the polynomial
$P$ may not be irreducible.

Some important problems on $q$-algebraic equations parallel those
in differential equations, and loosely speaking concern the existence, 
uniqueness, and some form of regularity of the solutions, all three problems 
being in terms of formal power series which may or may not be convergent. In
this setting, the regularity problem is to be interpreted as that of the
asymptotic behavior of the coefficients of the formal solutions. We will 
address these issues and emphazise on examples involving
divergent power series solutions.

Before putting our results in context, their flavor is illustrated on the
equation
$$\eqalignno{\qquad
  2f(z)=
  &{}4z^3f(qz)f(q^6z)+5q^2 z^6f(z)f(q^9z)f(q^{10}z)\cr
  &{}+18q^4 z^7 f(q^{14}z)^2+9z^{10}f(q^{-3}z)f(q^5z)f(q^{14}z)f(q^{16}z)\cr 
  &{}+3z^{14}f(q^{-5}z)+(q^8+2q^{17}) z^{14}f(z)\cr
  &{}+72z^{14} f(z)f(q^3z)f(q^5z)+1+15z^7 \, .
  &\equa{runningEx}\cr}
$$
From one of our result and
few lines of simple calculations it follows that whenever $q$ is a complex
number of modulus greater than $1$,
for some positive constant $c_q$ and for the positive root $R_q$ of the
equation $z^3+9q^{-24}z^7=1$,
$$
  f_n\sim c_q q^{n(n-3)} R_q^{-n}
  \eqno{\equa{runningExFinal}}
$$
as $n$ tends to infinity. The main interest of such formula is that $c_q$ does
not depend on $n$, so that the asymptotic growth of the sequence $(f_n)$ is
entirely captured by the term  $q^{n(n-3)} R_q^{-n}$. 

The significance of our results is related to the many applications 
of $q$-functional equations, and a section of the paper is devoted to examples
coming from combinatorics, dynamical systems, knot theory 
and statistical physics.

One of our results on the growth of $(f_n)$ is reminiscent of 
Maillet's (1903) theorem on formal power series
solutions of algebraic differential equations. Maillet's result was extended
to analytic nonlinear differential equations by Malgrange (1989), and though
our technique allows us to extend our result to some $q$-analytic functional 
equations, the needed assumptions on the analytic function are not particularly
sightly and we choose to restrict our exposition to $q$-algebraic equations,
a setting which seems suitable in most applications.

The Maillet-Malgrange theorem was extended to $q$-analytic functional
equations by Zhang (1998), who showed that the solution belong to {\it some}
$q$-Gevrey class, meaning that the sequence $(f_n)$ can grow at most as some
$q^{cn^2}$ for some unspecified $c$. This result was made more precise by
Di Vizio (2009) and Cano and Fortuny Ayuso (2012), in the sense that an upper
bound for $c$ was provided in terms of the Newton polygon associated with
a linearized version of the $q$-functional equation. While these results are
important from a theoretical perspective since they are indeed comparable to 
the Maillet-Malgrange theorem, they are insufficient for many applications, in 
particular in combinatorics and statistical physics where $(f_n)$ is 
the cardinality of some interesting
set, and one is interested in calculating this cadinality in a nice closed 
form. Obtaining such an explicit form is usually impossible with
our current knowledge, and one has to settle for less, such as
having a good estimate on the growth of $f_n$ with $n$. Our result
provides such estimate when $f$ is a divergent series, for in many cases 
it provides an explicit
asymptotic equivalent for $f_n$ up to a multiplicative constant.

A different view on the existence and regularity questions was developed
by Li and Zhang (2011) who gave sufficient conditions for an analytic
nonlinear $q$-difference equation to have an analytic solution. In the case
of $q$-algebraic equations, we provide an alternative sufficient condition.

It is worth noting that
the study of {\it linear} $q$-difference equations has a long history
spanning a century, from Carmichael (1912), Birkhoff (1913), 
Adams (1929, 1931), Trjitzinsky (1938) to the more recent works of 
B\'ezivin (1992), Ramis (1992), Sauloy (2000, 2003), Zhang (2002), Ramis, 
Sauloy and Zhang (2013). Nonlinear $q$-functional equations, which
seem more prevalent in applied mathematics, have been 
studied systematically only in the past few years by the mentioned authors.

Differing from the approach taken in the previous
papers, our results do not rely on the Newton polygon method, and the 
linearization of the equation that we will use is twofold: a partial one
on the equation, and one developed from the recursion
that the coefficients of the solution must satisfy. However, the full
potential of our method is achieved only when combined with the result
of Cano and Fortuny Ayuso (2012). Indeed, using their result, any $q$-algebraic
equation which has a solution in the Hahn field of generalized power series
can be transformed into one whose solution is a power series, and for which
the result of this paper may then be applied.

\bigskip

This paper is organized as follows. Section 2 contains the main theoretical
results. Section 3 discusses the computational aspects to implement our
results. Section 4 is devoted to some examples. Section 5 contains the proofs 
of the results stated in section 2. An appendix contains the {\tt Maple} code
which we used in one of the examples in section 4. The second section is 
divided in various subsections, the
first one introducing the mathematical objects which we will use, the second
one defining and discussing the notion of reduced equation, the third one
studying the existence of solutions of a reduced equation and stating some
basic regularity results. The fourth subsection provides an estimate of 
the growth of the coefficients for divergent solutions. The paper is 
essentially self-contained.

\bigskip

Throughout the paper, when $(u_n)$ is some sequence which does not vanish
ultimately and $c$ is a complex number, the relation $f_n\sim c u_n$, 
means that $\limn f_n/u_n=c$. In particular, if $c$ vanishes, this 
means $f_n=o(u_n)$ as $n$ tends to infinity.

If $f(z)=\sum_{n\geq 0} f_n z^n$ is a formal power series, we write $[z^n]f(z)$
or $[z^n]f$ for the coefficient $f_n$ of $z^n$ in this power series. The map
$[z^n]$ is linear on formal power series.

\bigskip

\noindent{\it Throughout this paper, $q$ is a complex 
number of modulus greater than $1$.}

\bigskip


\def\prevs{\the\sectionnumber.\the\subsectionnumber.\the\snumber }
\def\preveq{(\the\sectionnumber.\the\subsectionnumber.\the\equanumber)}

\section{Main results}
In order to state our results, the usual notation \qAlgebraic\ for
q-algebraic equation is not suitable. We will need a different formalism
which requires some definitions and notations which we will introduce in the
first subsection. In the second subsection, we will describe a way to bring
an algebraic $q$-functional equation into a suitable form. Our main
results are in the third and fourth subsections.

\subsection{Preliminaries}
Refering to the polynomial involved in \qAlgebraic, our most basic object will
be the monomials $z_0^{n_0}z_1^{n_1}\cdots z_k^{n_k}$, but we need to think of 
them in a different way, closer to their contribution in 
equation \qAlgebraic. As we will see afterwards, the following definition 
gives our analogue of the monomials.

\Definition{\label{qFactor}
  (i) A $q$-factor $A$ is a tuple $(a;\alpha_1,\ldots,\alpha_\ell)$ 
  in $\NN\times \ZZ^\ell$
  with $\alpha_1\leq\cdots\leq \alpha_\ell$ and $\ell$ is a positive integer.

  \medskip

  \noindent (ii) This $q$-factor is 
  shifting if $a$ is positive. A set of $q$-factors is shifting if all its
  elements are shifting.
}

\bigskip

When writing $q$-factors, a capital letter denotes the $q$-factor, the
corresponding lower case the first component of the $q$-factor, and the 
corresponding Greek one the other components. The index $\ell$ depends on
the $q$-factor. For instance, if $\calQ$ is a set of $q$-factors, a series
$\sum_{A\in\calQ} q^a z^{\alpha_\ell}$ is the sum over all $A$ in $\calQ$ of
$q$ at the power the first component of $A$ and $z$ at the power the last
component of $A$.

Definition \qFactor\ is introduced since to a $q$-factor $A$  corresponds 
an operator on formal power series
$$
  A f(z)=z^af(q^{\alpha_1}z)\cdots f(q^{\alpha_\ell}z) \, .
$$
We see that when $A$ is shifting,
$$
  [z^n]Af(z)=[z^{n-a}]f(q^{\alpha_1}z)\cdots f(q^{\alpha_\ell}z) 
$$
involves only $[z^i]f$ with $i\leq n-a\leq n-1$, hence the term `shifting'.

The set of $q$-factors can be viewed as a free $\CC$-module, and to a finite
weighted sum of $q$-factors, say $\sum r_i A_i$ corresponds the operator
$f\mapsto \sum r_i A_i f$.

Let $\calQ$ be a set of $q$-factors. In this paper we are 
primarily interested in the formal power series solutions of equations of 
\qAlgebraic, but written in the form
$$
  P(z)+\sum_{A\in\calQ} r_A\, A f(z)=0
  \eqno{\equa{qEq}}
$$
where $P(z)$ is a polynomial in $z$, the $r_A$ are some complex numbers, none
of them being $0$, and $\calQ$ is a finite set of $q$-factors.

\bigskip

\noindent{\bf Example.} Consider equation \runningEx. The set of $q$-factors 
involved in that equation is
$$\displaylines{\quad
  \calQ=\{\, (0;0), (3;1,6), (6;0,9,10), (7;14,14),
  \hfill\cr\hfill
   (10;-3,5,14,16), (14;-5), (14;0), (14;0,3,5) \,\}\, .
   \quad\equa{exampleQ}\cr}
$$
This $\calQ$ is not shifting since it contains the 
nonshifting $q$-factor $(0;0)$.

\bigskip

In light of the interpretation of $q$-factors as operators, the following
definition simply distinguishes the linear ones 

\Definition{\label{linearFactor}
  A $q$-factor $(a;\alpha_1,\ldots,\alpha_\ell)$ is linear if $\ell$ is $1$.
}

\bigskip

The nonshifting $q$-factors play a special role, which leads us to 
introduce the following notation.

\Notation{\label{QNotation}
  Let $\calQ$ be a set of $q$-factors. We write
  $$
    \calQ_0=\{\, A\in\calQ\,:\, a=0\,\} \, ,
  $$
  its subset of nonshifting $q$-factors and
  $$
    \calQ_+=\{\, A\in \calQ\,:\, a>0\,\}
  $$
  its subset of shifting ones.

  We set $\alpha(\calQ)=\max\{\, \alpha_\ell\,:\, A\in\calQ\,\}$.
}

\bigskip

Consider a $q$-algebraic equation \qEq\ without nonshifting $q$-factor. If 
$P(0)\not=0$, the equation has no solution. If $P(0)=0$, then we may 
factor $z$ out of the equation and simplify it by $z$. Therefore, for an
equation to have a solution, one must be able to simplify it into one that
has a nonshifting factor. Hence, there is no loss of generality in what follow
to consider only $q$-algebraic equations with at least one nonshifting factor. 

Considering the constant term
in \qEq, we see that for $f$ to be a solution we must have
$$
  P_0+\sum_{A\in\calQ_0} r_A f_0^\ell=0\, .
  \eqno{\equa{existenceCondition}}
$$
The only way for this equation not to have a solution is for the polynomial
$\sum_{A\in\calQ_0}r_A f_0^\ell$ to be degenerate and different than $-P_0$.
Equation \existenceCondition\ makes clear that several solutions may exist 
whenever $\calQ$ contains at least one nonshifting $q$-factor for 
which $\ell$ is at least $2$.


\subsection{Reduced  equations}
The equations in which all the nonshifting $q$-factors are linear are
important, and this motivates the following definition.

\Definition{\label{reduced}
  A $q$-algebraic equation is reduced if all its nonshifting $q$-factors are 
  linear.
}

\bigskip

Of course, not every $q$-algebraic equation is reduced. However, most of them
can be reduced as we will now explain. We will use an algorithm similar to
that proposed in Cano and Fortuny Ayuso (2012) to bring an equation to a 
quasi-solved form. However, since we are dealing exclusively with power
series, our algorithm is somewhat easier to describe.

We will write $\calO(z^k)$ to indicate a series in the ideal generated by
$z^k$. Thus, $f(z)=\calO(z^k)$ means that there exists a formal power
series $g(z)$ such that $f(z)=z^kg(z)$.

Assume that \existenceCondition\ has at least one solution $f_0$. We make the
change of function $f(z)=f_0+zg(z)$ in \qEq. If $A$ is nonshifting, then
$$\eqalign{
  Af(z)
  &{}=\prod_{1\leq i\leq \ell} \bigl( f_0+q^{\alpha_i}zg(q^{\alpha_i}z)\bigr) \cr
  &{}=f_0^\ell +f_0^{\ell-1} z\sum_{1\leq i\leq\ell} q^{\alpha_i} g(q^{\alpha_i}z)
    +\calO(z^2) \, . \cr
}
$$
If $A$ is shifting, then
$$
  Af(z)=z^a f_0^\ell+\calO(z^{a+1}) \, .
$$
Thus,
$$\displaylines{\qquad
  P(z)+\sum_{A\in\calQ} r_A Af(z)
  \hfill\cr\hfill
    {}=P_0+\sum_{A\in\calQ_0} r_A f_0^\ell 
    + z\sum_{A\in\calQ_0} r_A f_0^{\ell-1} \sum_{1\leq i\leq \ell} 
    q^{\alpha_i} g(q^{\alpha_i} z)
  \hfill\cr\hfill
  {}+\sum_{A\in\calQ_+} r_A z^a f_0^\ell + P_1z+\calO(z^2) \, .
  \qquad\equa{reduceA}\cr
  }
$$
Taking into account \existenceCondition, we see that the constant term
in \reduceA\ vanishes. Thus, after simplifying by $z$, we can rewrite \qEq\ as
$$
  \sum_{A\in\calQ_0} r_A f_0^{\ell-1} \sum_{1\leq i\leq \ell} 
  q^{\alpha_i}g(q^{\alpha_i}z)
  + \sum_{\ss A\in\calQ_+\atop\ss a=1} r_A f_0^\ell + P_1
  +\calO(z) =0 \, .
$$
The nonshifting part of this new equation is given by the operator
$$
  \sum_{A\in\calQ_0} r_A f_0^{\ell-1} \sum_{1\leq i\leq \ell} 
    q^{\alpha_i}(0;\alpha_i) \, .
$$
It involves only linear $q$-factors; however, this nonshifting 
operator may be $0$. In particular this is the case
if $f_0=0$ and $\ell>1$ for all nonshifting $q$-factors in the original 
equation. In this case, one may iterate this change of function, 
setting $g(z)=g_0+zh(z)$, and iterate further until one obtains a reduced 
equation. It seems that most equations encountered in applications can be
reduced in one or two changes of function. It is probably possible to 
characterize all equations that cannot be reduced by this change of function.

Note that for a reduced $q$-algebraic equation, condition \existenceCondition\
becomes $P_0+\sum_{A\in\calQ_0} r_A f_0=0$, so that a necessary condition
for such an equation to have a solution is that
$$
  P_0=0 \qquad\hbox{ or }\qquad \sum_{A\in\calQ_0} r_A\not=0 \, .
  \eqno{\equa{existenceConditionB}}
$$
If condition \existenceConditionB\ holds, then $f_0$ is unique 
if $\sum_{A\in\calQ_0} r_A\not= 0$ and is arbitrary otherwise.

\bigskip

\noindent{\bf Example.}{(qP${}_{\sevenrm I}$)} The $q$-difference first 
Painlev\'e equation, qP${}_{\hbox{\sevenrm I}}$,
is
$$
  \omega(qz)\omega\Bigl({z\over q}\Bigr)
  ={1\over \omega(z)}-{1\over z\omega^2(z)} \, .
$$
This equation, or an equivalent after a change of function, was introduced in 
Ramani and Grammaticos (1996) and further studied in the context of 
classification of rational surfaces by Sakai (2001), and from a different
perspective by Nishioka (2010) and Joshi (2012). The notation $\overline\omega$
and $\underline\omega$ for $\omega(qz)$ and $\omega(z/q)$ is convenient for
this example. This equation can be rewritten as
$$
  \underline \omega \omega^2\overline\omega=\omega-{1\over z} \, .
  \eqno{\equa{qPA}}
$$
This suggests to set $f(z)=\omega(1/z)$ to transform \qPA\ into
$$
  \underline f f^2\overline f=f-z \, ,
  \eqno{\equa{qPB}}
$$
which is a $q$-algebraic equation. The nonshifting $q$-factors are
$(0; -1,0,0,1)$ and $(0;0)$, and one of them, corresponding to the 
term $\underline f f^2\overline f$ is nonlinear. Applying $[z^0]$ to both
sides of \qPB, we must have $f_0^4=f_0$. Therefore, $f_0$ may be
$0$, $1$, $e^{2\imath \pi/3}$ or $e^{4\imath \pi/3}$. Following the indicated 
procedure to reduce the equation, we set $f(z)=f_0+zg(z)$.

\noindent{\it Case $f_0=0$.} In this case, $f(z)=zg(z)$ and \qPB\ becomes
$$
  z^3 \underline g g^2\overline g=g-1 \, .
  \eqno{\equa{qPC}}
$$
The only nonshifting $q$-factor is $(0;0)$, which corresponds to the linear
term $g$. Thus, this equation is reduced.

\noindent{\it Case $f_0=1$.} We then have $f(z)=1+zg(z)$ and \qPB\ becomes,
after some calculation
$$\displaylines{\qquad
 {1\over q} \ug + q\og +g +z\ug\og + {2\over q} z\ug g +2qz g\og +zg^2
  +2z^2\ug g \og 
  \hfill\cr\hfill
  {}+ {1\over q} z^2 \ug g^2 + qz^2 g^2\og
  +z^3\ug g^2\og +1 = 0 \, .
  \qquad\equa{qPD}\cr}
$$
The nonshifting $q$-factors are $(0;-1)$, $(0;1)$ and $(0;0)$, which are
all linear. Thus, this equation is reduced.

\noindent{\it Case $f_0=e^{2\imath \pi/3}$ or $f_0=e^{4\imath\pi/3}$.} 
Set $f(z)=f_0 h(z/f_0)$. One easily checks that $h$ solves \qPB\ and that
$h_0=1$, bringing us back to the previous case.

\bigskip

Following the custom in $q$-algebraic equations, we consider the operator
$\sigma$ and its powers $\sigma^n$ defined by
$$
  \sigma^n f(z)=f(q^n z) \, , \qquad n\in\ZZ \, .
$$

If $A$ is a $q$-factor, then
$$
  A\sigma^n f(z)=z^a f(q^{\alpha_1+n}z)\cdots f(q^{\alpha_\ell+n}z) \, ,
$$
so that $\sigma^n$ acts on $q$-factors to the right as
$$
  (a;\alpha_1,\ldots,\alpha_\ell)\sigma^n
  =(a;\alpha_1+n,\ldots,\alpha_\ell+n) \, .
$$
Consequently, if $f$ solves \qEq, then $g(z)=\sigma^{-n}f(z)$ solves
$$
  \sum_{A\in\calQ} r_A A\sigma^n g(z)+P(z) = 0 \, .
$$
Therefore, taking $n=-\alpha(\calQ_0)$, there is no loss of generality at all 
in assuming that $\alpha(\calQ_0)=0$.

\subsection{Existence and convergence or divergence of the solutions}
Not every reduced $q$-algebraic equation has a power series solution. For 
instance $zf(z)=1$ has not. Our first theorem is a sufficient condition for 
a reduced equation to have a solution. Recall that throughout the paper
we assume that $|q|>1$.

\Theorem{\label{existence}
  A sufficient condition for a reduced $q$-algebraic equation to have a 
  power series solution is that for any $n\geq 0$,
  $$
    \sum_{A\in\calQ_0} r_A q^{\alpha_1 n} \not =0 \, .
    \eqno{\equa{existenceConditionA}}
  $$
  If this condition holds, the solution is uniquely determined by $f_0$.
}

\bigskip

Condition \existenceConditionA\ is analoguous to the condition introduced
by Cano and Fortuny Ayuso (2012) for an equation in quasi-solved form to be
in solved form, and the second 
assertion of Theorem \existence\ may be seen as a rephrasing of their Theorem
2 in our context of power series.

\bigskip

\noindent{\bf Example.} For equation \runningEx, which is reduced, we have
$$
  \sum_{A\in\calQ_0} r_A q^{\alpha_1 n}
  = 2
$$
and condition \existenceConditionA\ holds.

For the $q$-difference first Painlev\'e equations \qPC,
the unique nonshifting power is $(0;0)$, which corresponds to the term $g$
in that equation, so that \existenceConditionA\ is $1\not=0$.

For equation \qPD, condition \existenceConditionA\ is
$$
  {1\over q^{n+1}}+q^{n+1}+1\not= 0\, .
$$
This condition is always satisfied when $|q|>1$ since the polynomial
$z^2+z+1$ has two roots of modulus $1$.

\bigskip

Once a reduced equation has a solution under 
criterion \existenceConditionA, the question arises as to whether this solution
has a positive radius of convergence or is a divergent series.

If $\calQ_+$ is empty, a reduced equation has the form
$$
  \sum_{A\in\calQ_0} r_A f(q^{\alpha_1}z) +P(z)=0 \, .
$$
Applying $[z^n]$ to both sides of this equation,
$$
  \Bigl( \sum_{A\in\calQ_0} r_A q^{\alpha_1 n} \Bigr) f_n+P_n=0 \, .
$$
Under condition \existenceConditionA, we obtain $f_n$ and since $p$ is the
degree of $P$,
$$
  f(z)=-\sum_{0\leq n\leq p} {P_n\over \sum_{A\in\calQ_0} r_A q^{\alpha_1 n}}
  z^n
$$
is a polynomial in $z$. In what follows we will then assume that $\calQ_+$
is not empty. Our next result provides a simple test for convergence of 
the solutions. To avoid some trivial situations, we use the following 
definition.

\Definition{\label{collected}
  We say that the $q$-factors involved in a $q$-algebraic equation have been
  collected if each $q$-factor occurs only once in the equation.
}

\bigskip

For instance, the $q$-factors $(0;1)$ and $(0;0)$ in the
equation $f(qz)+f(z)-f(z)=0$ have not been collected. There is of course no
loss of generality to assume that the $q$-factors have been collected.


\Theorem{\label{convergenceDivergence}
  Consider a reduced $q$-algebraic equation whose $q$-factors have been
  collected, such that \existenceConditionA\ 
  holds and $\calQ_+\not=\emptyset$. If $\alpha(\calQ_0)\geq \alpha(\calQ_+)$
  then the unique power series solution has a positive radius of convergence.
}

\bigskip

It is possible for the radius of convergence to be infinite, as in the
equation
$$
  f(z)+qzf(z)-z-qz^2=0
$$
which is designed to have $f(z)=z$ for solution.

In light of Theorem \convergenceDivergence, note that a $q$-algebraic 
equation may have both convergent and divergent solutions, because different
initial values $f_0$ usually lead to different reduced equations.

\bigskip

\noindent{\bf Example.} {(qP${}_{\sevenrm I}$)} The $q$-factors involved in \qPD\ are
$$\displaylines{\quad
  \calQ=\{\, (0;-1), (0;0), (0;1), (1;-1,1), (1; -1,0); (1; 0,1), 
  \hfill\cr\hfill
  (1;0,0), (2;-1,0,1), (2; -1, 0, 0), (2; 0,0,1); (3; -1,0,0,1)\,\} \, .
  \cr}
$$
Thus $\alpha(\calQ_0)=1=\alpha(\calQ_+)$ and Theorem \convergenceDivergence\
implies that the solution has a positive radius of convergence.

\bigskip

Theorem \convergenceDivergence\ raises the question of what can be said
when $\alpha(\calQ_0)<\alpha(\calQ_+)$. This question is partially answered
in the next subsection.

\subsection{Asymptotic behavior of the coefficients of the divergent solutions}
In this subsection we will derive an asymptotic equivalent for
the coefficients of the solution of a reduced $q$-algebraic equation
satisfying \existenceConditionA. To state our results, we need to introduce 
further definitions related to $q$-factors.

\Definition{\label{qHeight}
  Consider a set $\calQ$ of $q$-factors such that $\calQ_+\not=\emptyset$ and
  $\alpha(\calQ_0)=0$.

  \medskip

  \noindent(i) The height of a shifting $q$-factor 
  $A=(a;\alpha_1,\ldots,\alpha_\ell)$ is
  $$
    H(A)=\alpha_\ell/(2a) \, .
  $$

  \noindent
  (ii) If $\calQ$ is a finite set of $q$-factors, its height is the 
  largest height of its shifting elements, that is
  $$
    H(\calQ)=\max_{A\in\calQ_+} H(A) \, .
  $$

  \noindent
  (iii) The crest $\widehat\calQ$ of $\calQ$ is
  $$
    \widehat\calQ
    =\{\,  A\in\calQ\,:\, 2aH(\calQ)=\alpha_\ell\,\}\, .
  $$

  \noindent
  (iv) The co-height of $\calQ$ is
  $$
    h(\calQ)=\min\{\, a\,: A\in\widehat\calQ\cap \calQ_+\,\} \, .
  $$

  \noindent
  (v) The scope of $A$ is the number of maximal $\alpha_i$,
  $$
    s(A)=\sharp\{\, i\,:\, \alpha_i=\alpha_\ell\,\} \, .
  $$

  \noindent
  (vi) The $\calQ$-Borel transform of a formal power 
  series $f(z)=\sum_{n\geq 0} f_n z^n$ is the formal power series
  $$
    \calB_\calQ f(z)=\sum_{n\geq 0} q^{-H(\calQ) n(n-h(\calQ))} f_n z^n \, .
  $$
}

The crest $\widehat\calQ$ coincides with the set of shifting $q$-factors of 
maximal height together with the set of nonshifting $q$-factors for which 
$\alpha_\ell$ vanishes. Given how the crest is defined, we see 
that $2aH(\calQ)-\alpha_\ell$ is positive for all $q$-factors not 
in $\widehat\calQ$.

\bigskip

\noindent{\bf Example.} The heights of the shifting $q$-factors in \exampleQ\ 
are $1$, $5/6$, $1$, $8/10$, $-5/28$, $0$, $5/28$. The largest height is
$H(\calQ)=1$. Since both $(3;1,6)$ and $(7;14,14)$ have height $1$, the 
crest of $\calQ$ is then $\widehat\calQ=\{\, (0;0), (3;1,6), (7;14,14)\,\}$ . 
The co-height of $\calQ$ is $h(\calQ)=\min(3,7)=3$.

\bigskip

To each reduced $q$-functional equation of the form \qEq\ for 
which $\calQ_+$ is not empty, we associate a polynomial as follows.

\Definition{\label{crestPoly}
  Assume that \qEq\ is reduced and that $\alpha(\calQ_0)=0$.
  The crest polynomial associated to \qEq\ is
  $$
    \calC_{q,t}(z)
    =\sum_{A\in\widehat\calQ} r_A s(A) q^{-H(\calQ)a(a-h(\calQ))}z^a t^{\ell-1}\, .
  $$
}

There is a slight abuse of terminology in Definition \crestPoly\ for
a $q$-algebraic equation can be divided by any nonzero complex number. So
the crest polynomial is defined only up to a multiplicative constant.
However, since we will be interested in the zeros of this polynomial, this
abuse of terminology will not create any ambiguity.

We see that the crest polynomial is the $\calQ$-Borel transform of
$\sum_{A\in\widehat\calQ} r_A s(A)z^a t^{\ell-1}$.

From the Definition \qHeight.iii, if $A$ is in $\widehat\calQ$ and $a=0$, 
then $\alpha_\ell=0$, and $\ell=1$ if the equation is reduced. Therefore, the 
constant term of the crest polynomial is
$$
  \calC_{q,t}(0)=\sum_{A\in\calQ_0\cap\widehat\calQ} r_A s(A) \, .
$$
If \qEq\ is reduced, then $\calQ_0$ contains only linear $q$-factors, which
forces $s(A)=1$ whenever $A$ in nonshifting, and then
$$
  \calC_{q,t}(0)=\sum_{\ss A\in\calQ_0\atop\ss\alpha_\ell=0} r_A \, .
$$
Thus, condition \existenceConditionA\ for $n=0$ may be
rewritten as $\calC_{q,t}(0)\not=0$.

\bigskip
\noindent{\bf Example.} To evaluate the crest polynomial associated
to \runningEx, we have the following table listing the parts of the equation
relevant to the crest.
\bigskip

\setbox1=\vbox{\halign{$#$\hfill\quad&\hfill$#$\quad&\hfill$#$\quad&\hfill$#$\cr
A & (0;0) & (3;1,6) & (7;14,14) \cr
r_A& -2  & 4        & 18q^4 \cr
s(A) & 1 & 1   & 2\cr
\ell & 1 & 2 & 2 \cr
H(\calQ)a\bigl(a-h(\calQ)\bigr) & 0 & 0 & 28 \cr}}

\centerline{\box1}
\bigskip
\noindent
The crest polynomial is then
$$
  \calC_{q,t}(z)=-2+4z^3t+36q^{-24}z^7t \, .
$$


If $A$ is in the crest of $\calQ$, its height is that of $\calQ$, so that
$2H(\calQ)a\bigl(a-h(\calQ)\bigr)=\alpha_\ell \bigl(a-h(\calQ)\bigr)$ is an
integer. Also, $\ell$ is at least $1$. Therefore, $\calC_{q,t}(z)$ is a 
polynomial in $(z,q^{-1/2},t)$, provided we ignore a possible dependence
in $q$ in the coefficient $r_A$.

If $\calQ$ is shifting, the
polynomial $\sum_{A\in\widehat\calQ\setminus\calQ_0} r_A q^{-H(\calQ)a(a-h(\calQ))}z^a$
is in the ideal generated by $z$. However, it is possible that this sum
vanishes due to some cancelations. Consequently, a crest polynomial
may be constant. If this is the case, we consider that it has a root at 
infinity for the following result to hold without further discussion.

Our main result is the following, and its implications will be discussed
afterwards.

\Theorem{\label{mainTh}
  Consider a reduced $q$-algebraic equation where the $q$-factors have been
  collected, such that \existenceConditionA\
  holds and $0=\alpha(\calQ_0)<\alpha(\calQ_+)$. Let $f$ be a solution of \qEq. 
  Let $\calC_{q,t}(z)$ be the crest polynomial 
  associated to \qEq, and let $R_{q,f_0}$ be the smallest modulus of the zeros
  of $\calC_{q,f_0}$.

  \medskip
  \noindent (i) $\calB_\calQ f(z)$ has radius of convergence $R_{q,f_0}$.

  \medskip
  \noindent (ii) There exists some positive $\Theta$ such that
  $\calC_{q,f_0}(z)\calB_\calQ f(z)$ is a convergent power series that has no 
  singularities other than removable ones in the disk centered at $0$ and of 
  radius $q^\Theta R_{q,f_0}$.

  \medskip
  \noindent (iii) if
  $$
    \{\, P_i\,:\, 0\leq i\leq p\,\}\cup\{\, r_A\,:\, A\in\calQ_+\,\}
    \cup\{\, -r_A\,:\, A\in\calQ_0\,\}
  $$
  is a set of real numbers all of the same sign, then $(f_n)$ is a nonnegative 
  sequence. Moreover, if $\calQ_0=\{\, (0;0)\,\}$ then 
  $\calC_{q,f_0}(z)\calB_\calQ f(z)$ is a power series whose 
  coefficients are nonnegative provided we choose $r_{(0;0)}$ nonnegative, 
  and these coefficient all vanish if and only if $f=0$.
}

\bigskip

In the statement of Theorem \mainTh, recall that $\alpha(\calQ_0)=0$
is not restrictive at all, and that Theorem \convergenceDivergence\ asserts 
that $\alpha(\calQ_+)>0$ is the only case where we can have a divergent 
solution. The third assertion of Theorem \mainTh\ implies that the solution
is indeed divergent in equations that have coefficients of the proper signs
and a single nonshifting $q$-factor. Recall that one should 
read $R_{q,f_0}=+\infty$ if $\calC_{q,f_0}$ is a nonzero constant polynomial 
and $R_{q,f_0}=0$ if $\calC_{q,f_0}$ is the constant polynomial $0$.

\bigskip

Our proof shows that $\Theta$ in statement (ii) of Theorem \mainTh\ may be 
taken to be the smallest of $2H(\calQ)$,
$\max_{A\in\calQ\setminus\widehat\calQ}(2Ha-\alpha_\ell)$ and $1$.

A careful examination of the proof shows that the power series 
$(\calC_{q,f_0}\calB_\calQ f)(z)$ is almost a linear combination of tangled 
products in the sense of Garsia (1981) 
of $\calQ$-Borel transforms of $f$. However, since there is no simple 
way to calculate this $\calQ$-Borel transform, such expression does 
not appear useful to estimate $[z^n]f$.

The strength of Theorem \mainTh\ comes from the following. Set 
$$
  U_{q,f_0}(z)= \calC_{q,f_0}(z)\calB_\calQ f(z) \, .
$$
Theorem \mainTh\ asserts that $U_{q,f_0}$ has no singularity of modulus 
less than $q^\Theta R_{q,f_0}$. The coefficient
$$
  q^{-H(\calQ)n(n-h(\calQ))} f_n
  = [z^n]{U_{q,f_0}(z)\over \calC_{q,f_0}(z)}
  \eqno{\equa{basicFormula}}
$$
can be evaluated via the Cauchy 
formula, and results on the asymptotic behavior of these coefficients are
readily available through the singularity analysis of meromorphic 
functions, as explained in Flajolet and Sedgewick (2009; see in particular 
Chapter IV).

\bigskip

\noindent{\bf Example.} For equation \runningEx, $f_0=1/2$ and 
the crest polynomial is
$\calC_{q,1/2}(z)=-2+2z^3+18q^{-24}z^7$ has a positive solution, call it $R$,
which then satisfies
$$
  1=R^3+9q^{-24}R^7 \, .
$$
To see that $R$ is the solution of smallest modulus, if $\calC_{q,1/2}(Rz)=0$ 
had another solution $\zeta$ of modulus at most $1$, then
we would have
$$
  1=R^3\zeta^3+9q^{-24}R^7\zeta^7
$$
so that we would be able to write $1$ as a convex combination of $\zeta^3$
and $\zeta^7$ which both belong to the convex unit disk centered at $0$ and
of radius $1$ and have relatively prime exponents; this forces $\zeta=1$.

We then write $\calC_{q,1/2}(z)=(1-z/R)Q(z)$ where $Q$ is some polynomial of 
degree $6$. It then follows from singularity analysis (see Flajolet and 
Sedgewick, 2009; chapter IV) that, with the notation as in Theorem \mainTh,
$$
  [z^n]{U_{q,1/2}(z)\over\calC_{q,1/2}(z)} 
  = [z^n]{U_{q,1/2}(z)\over (1-z/R)Q(z)}
  \sim R^{-n} U_{q,1/2}(R)/Q(R)
  \eqno{\equa{exU}}
$$
as $n$ tends to infinity. Set $c_q=U_{q,1/2}(R)/Q(R)$. Theorem \mainTh.iii
implies that $U_{q,1/2}\not=0$, hence $c_q$ does not vanish. This 
yields \runningExFinal.

\bigskip

Following B\'ezivin (1992), Ramis (1992) and others, recall that a formal
power series $f(z)=\sum_{n\geq 0} f_n z^n$ is of $q$-Gevrey order $s$
if the series $\sum_{n\geq 0} f_n q^{-sn^2/2}z^n$ has a positive radius of
convergence. Since we consider $|q|>1$, if $f$ has $q$-Gevrey order $s$, 
it also has $q$-Gevrey
order any number greater than $s$. It is then of interest to find the smallest
$q$-Gevrey order of a divergent power series. Clearly,
Theorem \mainTh\ yields some information on the $q$-Gevrey
order of the solution of \qEq, which is sharp when the assumptions in
Theorem \mainTh.iii are satisfied.

\Corollary{\label{qGevrey}
  Under the assumptions of Theorem \mainTh, the power series solution of
  \qEq\ has $q$-Gevrey order at least $2H(\calQ)$. Moreover, if 
  assumptions of Theorem \mainTh.iii hold, then $2H(\calQ)$
  is the smallest $q$-Gevrey order of the power series solution.
}


\Remark The reduction algorithm described in section \fixedref{2} changes the
equation and, therefore, may change its height. In order for Theorem \mainTh\
and Corollary \qGevrey\ to deliver a sharp result, we should reduce the
equation in such a way that its height is as small as possible. Assume that
the solution $f(z)$ of equation \qEq\ is such that $f_0=0$. We then set
$f(z)=zg(z)$. For a $q$-factor $A=(a;\alpha_1,\ldots,\alpha_\ell)$,
$$
  Af(z)=q^{\alpha_1+\cdots+\alpha_\ell}z^\ell Ag(z) \, .
$$
We then set
$$
  \widetilde A=(a+\ell-1;\alpha_1,\ldots,\alpha_\ell)
$$
and $\tilde r_A=q^{\alpha_1+\cdots+\alpha_\ell}r_A$, so that $g$ is a solution of
$$
  P(z)+\sum_{A\in\calQ}\tilde r_A z \widetilde A g(z)=0 \, .
$$
Since $f_0=0$, we have $P(0)=0$ and we can simplify this equation by $z$
to obtain
$$
  {P(z)\over z}+\sum_{A\in\calQ}\tilde r_A \widetilde A g(z)=0 \, .
$$
We have
$$
  H(\widetilde A)={\alpha_\ell\over 2(a+\ell-1)}\leq H(A)
$$
and the inequality is strict whenever $\ell>1$. Thus, when $f_0=0$, we
should reduce the equation, in particular if the crest of $\calQ$ contains
nonlinear $q$-factors, that is factors for which $\ell>1$.

\bigskip

The following corollary to Theorem \mainTh\ covers important applications 
and shows that some oscilatory behavior may occur and gives a sharp bound
on the periodicity.

\Corollary{\label{mainCor}
  Under the assumptions of Theorem \mainTh, if the crest $\widehat\calQ$ has a 
  unique shifting element $A=(a;\alpha_1,\ldots,\alpha_\ell)$, then
  there exist some complex numbers $c_0,\ldots,c_{a-1}$ such that, as $n$ tends
  to infinity,
  $$
    [z^n]f\sim q^{H(\calQ)n(n-h(\calQ))} 
    \Bigl({-r_A s(A)f_0^{\ell-1}\over r_{(0;0)}}\Bigr)^{n/h(\calQ)} c_m
  $$
  where $m$ is the remainder in the Euclidean division algorithm of $n$ by $a$.
}


\Remark In the statement of Corollary \mainCor, $r_As(A)f_0^{\ell-1}/r_{(0;0)}$
may not be a positive real number. Thus, to take the fractional 
power $n/h(\calQ)$ of this number requires one to choose a branch cut in the 
complex plane. The fractional power is then defined up to some 
$e^{2i\pi kn/h(\calQ)}$ for some $0\leq k<h(\calQ)$. Since $h(\calQ)=a$, this
indeterminacy may be absorbed in the constant $c_m$.

As far as the oscillatory behavior of $[z^n]f$ is concerned, 
Corollary \mainCor\ shows that it may be decomposed into two parts: write
$r_As(A)f_0^{\ell-1}/r_{(0,0)}$ as $\rho e^{2i\pi \theta}$ with $\rho$ nonnegative
and $\theta$ in $[\,0,1)$. Corollary \mainCor\ asserts that
$$
  [z^n]f\sim q^{H(\calQ)n(n-h(\calQ))} \rho^{n/h(\calQ)} 
  e^{2i\pi\theta n/h(\calQ)} c_m 
$$
as $n$ tends to infinity. If $\theta$ is irrational, the 
sequence $(e^{2i\pi\theta n/H(\calQ)})_{n\in\NN}$ is dense in the unit circle, 
making $[z^n]f$ with seemingly little regularity. If $\theta$ is rational, 
say $\theta=h(\calQ)p'/p$ with $p'$ and $p$ positive integers mutually 
prime, the 
sequence $(e^{2i\pi\theta n\theta/h(\calQ)})_{n\in\NN}$ has periodicity $p$. Thus,
the sequence $(e^{2i\pi n\theta/h(\calQ)}c_m)_{n\in\NN}$ has periodicity the 
least common multiple of $p$ and $a$. Example \fixedref{5}
of section \fixedref{5} will illustrate this phenomenon with $p=4$ and $a=17$,
leading to a periodicity of $68$. However, because $(c_m)_{0\leq m<a}$ may
have a periodicity a divisor of $a$, it is possible that the 
sequence $(e^{2i\pi n\theta/h(\calQ)}c_m)_{n\in\NN}$ has a period smaller than the
least common multiple of $p$ and $a$.

\bigskip

\Proof Since the crest has a unique shifting element $A$ 
and $\alpha(\calQ_0)=0$, the height of $\calQ$ is $\alpha_\ell/2a$
and its co-height is $a$. Since $\alpha(\calQ_0)=0$, the crest contains
exactly the elements $(0;0)$ and $A$. The crest polynomial is
$$
  \calC_{q,f_0}(z)=r_{(0;0)}+r_A s(A)z^a f_0^{\ell-1} \, .
$$
Its roots are some $\zeta_0$ and $\zeta_k=e^{2ik\pi/a}\zeta_0$, $1\leq k<a$.
They all have the same modulus $|r_{(0;0)}/r_A s(A)f_0^{\ell-1}|^{1/a}$. 
Through \basicFormula, Theorem \mainTh\ yields
$$
  [z^n]f
  ={q^{(\alpha_\ell/2a) n(n-a)}\over r_{(0;0)}} [z^n]
  { U_{q,f_0}(z)\over\prod_{0\leq k<a}(1-z/\zeta_k)} \, .
$$
But
$$\eqalign{
  {1\over \prod_{0\leq k<a}(1-z/\zeta_k)}
  &{}= \sum_{0\leq k<a}
    { 1\over 1-z/\zeta_k }
    \prod_{\ss 0\leq j<a\atop\ss j\not=k}{1\over 1-\zeta_k/\zeta_j}\cr
  &{}=\sum_{0\leq k<a}
    { 1\over 1-z/\zeta_k }
    \prod_{\ss 0\leq j<a\atop\ss j\not=k}{1\over 1-e^{2i(k-j)\pi/a}}
  \, .\cr}
$$
Thus,
$$
  [z^n]f\sim {q^{H(\calQ)n(n-h(\calQ))}\over r_{(0;0)}} \sum_{0\leq k<a} 
  {U_{q,f_0}(\zeta_k)\over \prod_{\ss 0\leq j<a\atop\ss j\not=k} 1-e^{2i(k-j)\pi/a}}
  \zeta_k^{-n} \, .
$$
Note that $(\zeta_k/\zeta_0)^{-n}$ depends only on the remainder $m$ in
the Euclidean division algorithm of $n$ by $a$. We then set
$$
  \tilde c_m={1\over r_{(0;0)}}\sum_{0\leq k<a}
  {U_{q,f_0}(\zeta_k)\over \prod_{\ss 0\leq j<a\atop\ss j\not=k} 1-e^{2i(k-j)\pi/a}}
  \Bigl({\zeta_0\over \zeta_k}\Bigr)^n \, ,
$$
and $c_m=e^{2\imath k\pi/a}\tilde c_m$ for some $0\leq k<a$ according to the 
choice of $\zeta_0$  among the $a$ roots of the crest polynomial.\hfill\qed

\bigskip

\noindent{\bf Example.} {\it (qP${}_I$, continued)} Consider the
$q$-functional equation \qPC. One can see by induction that $g$ is in fact
a function of $z^3$ (see also Joshi, 2012). Thus, we set $g(z)=h(z^3)$
and rewrite \qPC\ as
$$
  -h(z^3)+z^3h\Bigl({z^3\over q^3}\Bigr) h(z^3)^2 h(q^3z^3)+1=0 \, .
$$
Setting $r=q^3$, and substituting $z$ for $z^3$, we have
$$
  -h(z)+z h(z/r)h(z)^2 h(rz)+1=0 \, ,
$$
which, in the notation of $q$-factors means
$$
  \bigl( -(0;0)+(1;-1,0,0,1)\bigr) h(z) +1=0 \, .
$$
For this equation, $H(\calQ)=1/2$ and $h(\calQ)=1$. Corollary \mainCor\ yields
$$
  h_n\sim c q^{n(n-1)/2}
$$
as $n$ tends to infinity. It follows that $g_{3n+1}=g_{3n+2}=0$ and
$$
  g_{3n}=h_n\sim c q^{n(n-1)/2} \, .
$$
Going back to equation \qPB, we then have $f_n=g_{n-1}$ which yields
$f_{3n}=f_{3n+2}=0$ and $f_{3n+1}\sim c q^{n(n-1)/2}$. Hence
$\omega(z)=\sum_{n\geq 0} \omega_n z^{-n}$ with $\omega_{3n}=\omega_{3n+2}=0$
and $\omega_{3n+1}\sim c q^{n(n-1)/2}$ as $n$ tends to infinity.

\bigskip


\def\prevs{\the\sectionnumber.\the\snumber }
\def\preveq{(\the\sectionnumber.\the\equanumber)}

\section{Computational aspects} 
The reduction procedure often increases the number of terms in a
$q$-algebraic equation, as it may be seen for instance in comparing
\qPB\ with \qPD. In the next section, devoted to examples, we will
consider an equation due to Cano and Fortuny Ayuso (2012) for which
the reduction procedure needs to be applied 10 times, leading to an
equation with 397 terms. Such an example makes clear that for our
theory to be useful in applications, it has to be implemented on a
computer algebra package. The goal of this section is not to describe
an implementation for a specific computer algebra system, but to
present algorithms from which an efficient implementation is easily
written. Since it reformulates our theory in terms of multivariate
polynomials, it also gives it a more algebraic-geometric flavor.

To a $q$-factor $(a;\alpha_1,\ldots,\alpha_\ell)$ corresponds a
monomial $z^aY_{\alpha_1}\cdots Y_{\alpha_\ell}$. We will also consider
monomials as operators, with then 
$$
  (z^aY_{\alpha_1}\cdots Y_{\alpha_\ell})f(z)
  =z^af(q^{\alpha_1}z)\cdots f(q^{\alpha_\ell} z) \, .
$$
In particular, following one of the many known ways of writing $q$-algebraic
equations, \qAlgebraic\ may be rewritten as
$$
  P(z,Y_0,\ldots,Y_k)f(z)=0 \, .
  \eqno{\equa{algoA}}
$$
Recall that the total degree of a monomial is
$$
  \Deg (z^aY_0^{n_1}\cdots Y_k^{n_k})=a+n_1+\cdots+n_k \, ,
$$
and that the total degree of a polynomial is the largest degree of its
monomials. We can now transcribe the meaning of a reduced equation.

\Proposition{\label{algoReduced}
  Equation \algoA\ is reduced if and only if $\Deg P(0,Y_0,\ldots,Y_k)=1$.
}

\bigskip

\Proof Using the correspondence between monomials and $q$-factors,
$$
  [z^0]P(z,Y_0,\ldots,Y_k)
  = \sum_{\ss A\in\calQ\atop\ss a=0} r_A A + P(0,\ldots,0)
  = P(0,Y_0,\ldots,Y_k) \, .
$$
Thus, $P(0,Y_0,\ldots,Y_k)$ has total degree $1$ if and only if each
$A$ in $\calQ$ with $a=0$ is of the form $(0;\alpha_\ell)$, and therefore
is linear.\hfill\qed

\bigskip

In the reduction steps, we set $f(z)=f_0+zg(z)$. This corresponds to 
transforming the equation as follows.

\Proposition{\label{algoChangeFunction}
  The power series $f(z)=f_0+zg(z)$ solves \algoA\ if and only if
  $P(0,f_0,\ldots,f_0)=0$ and
  $$
    P(z,f_0+zq^0Y_0,f_0+zq^1Y_1,\ldots,f_0+zq^kY_k)g(z)=0 \, .
  $$
}

\Proof We have
$$
  [z^0]P\bigl(z,f(z),\ldots,f(q^kz)\bigr)=P(0,f_0,\ldots,f_0) \, .
$$
The second assertion follows from the fact that
$$
  Y_j\bigl( f_0+zg(z)\bigr)
  = f_0+zq^jg(q^jz)
  = f_0+zq^jY_jg(z) \, .
  \eqno{\qed}
$$

Propositions \algoReduced\ and \algoChangeFunction\ yield the following 
algorithm for reducing a $q$-algebraic equation. We first define a
procedure which removes trivial factors of a polynomial whenever possible.
Recall that the order of a polynomial with respect to a variable $Y_j$ is 
the largest $n$ such that $Y_j^n$ divides the polynomial. We write
$\ord_{Y_j}$ for that order.

\bigskip

\algorithm{%

Procedure RemoveTrivialFactors(polynomial $R$)

return $R/(z^{\ord_z R}Y_0^{\ord_{Y_0}R}\cdots Y_k^{\ord_{Y_k}R})$

}

\bigskip

\noindent and the algorithm is

\bigskip

\algorithm{%

$P(z,Y_0,\ldots,Y_k)\leftarrow$ RemoveTrivialFactors$\bigl(P(z,Y_0,\ldots,Y_k)\bigr)$

while $\Deg P(0,Y_0,\ldots,Y_k)>1$ do

\quad  solve for $f_0$ in $P(0,f_0,\ldots,f_0)=0$

\quad  $P(z,Y_1,\ldots,Y_k) \leftarrow$ RemoveTrivialFactors$\bigl(P(z,f_0+zq^0Y_0,$

\quad\hfill $f_0+zq^1Y_1,\ldots,f_0+zq^kY_k)\bigr)$\par

}

\bigskip

If this loops terminates, which it may not in rather unusual cases, then the 
equation $P(z,Y_0,\ldots,Y_k)f(z)=0$ is reduced. Thus, from now on we assume
that the equation is reduced.\note{can we write down when the loop does
not terminate?}

Note that the simplicity of the algorithm as it is written hides that
the equation $P(0,f_0,\ldots,f_0)$ may have multiple roots, so that if one
would like to keep track of all the solutions, an extra tree-like structure
is needed.

At the end of subsection \fixedref{2.2}, we mentioned the important aspect
that there is no loss of generality to assume $\alpha(\calQ_0)=0$, since
we can make a change of function $g(z)=\sigma^{-n}f(z)$. This property is
awkward on polynomials, because the indices of the variables $Y_j$ have
meaning, and can be rewritten as follows.

\Proposition{\label{algoqMult}
  The formal power series $f=\sigma^ng$ solves \algoA\ if and only if
  $$
    P(z,Y_n,\ldots,Y_{k+n})g(z)=0\, .
  $$
}

\Proof We have $Y_if(z)=f(q^iz)=g(q^{i+n}z)=Y_{i+n}g(z)$.\hfill\qed

\bigskip

We then see that we may need to use polynomials in variables with negative
indices, which is a substantial distinction with the usual algebraic
geometry where the names of variables are irrelevant. So instead of a 
polynomial
in the $k+2$ variables $z,Y_0,\ldots,Y_k$, we will consider one in $2k+2$
variables yielding the equation
$$
  P(z,Y_{-k},\ldots,Y_k)f(z)=0 \, ,
  \eqno{\equa{algoB}}
$$
and assume that this equation is reduced.

We set
$$
  P_0(Y_{-k},\ldots,Y_k)=P(0,Y_{-k},\ldots,Y_k)-P(0,\ldots,0)
$$
and
$$\displaylines{\qquad
  P_+(z,Y_{-k},\ldots,Y_k)
  =P(z,Y_{-k},\ldots,Y_k)-P_0(Y_{-k},\ldots,Y_k)
  \hfill\cr\hfill
  {}-P(z,0,\ldots,0) \, .\qquad\cr}
$$
The notation $P_0$ and $P_+$ indicates that these polynomials represent
the nonshifting and the shifting part of the equation since
$$
  \sum_{A\in\calQ_0} r_A Af(z)=P_0(Y_{-k},\ldots,Y_k)f(z)
$$
and
$$
  \sum_{A\in\calQ_+} r_A Af(z)=P_+(z,Y_{-k},\ldots,Y_k)f(z) \, .
$$
Then, $\alpha(\calQ_0)$ is the largest $j$ such that $P_0(Y_{-k},\ldots,Y_k)$ 
does not depend on $Y_{j+1},\ldots,Y_k$. It can be calculated with the
following obvious algorithm.

\bigskip

\algorithm{

$j\leftarrow k$

while ${\ds\partial\over\partial\ds Y_j} P_0(0,Y_{-k},\ldots,Y_k)=0$ do $j\leftarrow j-1$

$\alpha(Q_0)\leftarrow j$
}

\bigskip

Similarly, $\alpha(\calQ_+)$ is the largest $j$ such that
$P_+(Y_{-k},\ldots,Y_k)$ does not depend on $Y_{j+1},\ldots,Y_k$ and may
be calculated with the following algorithm.

\bigskip

\algorithm{

$j\leftarrow k$

while ${\ds\partial\over\ds\partial Y_j}P_+(z,Y_{-k},\ldots,Y_k)=0$ do $j\leftarrow j-1$

$\alpha(Q_+)\leftarrow j$
}

\bigskip

We can restate the analogue of condition \existenceCondition.

\Proposition{\label{algoExistenceCond}
  For the reduced equation \algoB, condition \existenceConditionA\ is equivalent
  to
  $$
    P(0,q^{-kn},\ldots,q^{kn})\not=0 \qquad\hbox{for any $n\geq 0$.}
    \eqno{\equa{algoExistenceConditionA}}
  $$
}

\Proof Since the equation is reduced,
$$\eqalign{
  \sum_{A\in\calQ_0} r_A f(q^{\alpha_\ell}z)
  &{}=\sum_{-k\leq i\leq k} \sum_{\ss A\in\calQ_0\atop\ss \alpha_\ell=i}
    r_A f(q^iz) \cr
  &{}=\sum_{-k\leq i\leq k} [Y_i]P_0(Y_{-k},\ldots,Y_k)f(q^iz) \cr
  &{}=P_0\bigl(f(q^{-k}z),\ldots,f(q^kz)\bigr) \, .\cr}
$$
The result follows by taking $f(z)=z^n$.\note{should we try to translate this into an algorithm? Seems complicated}\hfill\qed

\bigskip

Then we can restate Theorem \convergenceDivergence.


\Proposition{\label{algoCvDv}
  Consider a reduced equation \algoB\ such that \algoExistenceConditionA\ 
  holds and $\alpha(\calQ_0)\geq \alpha(\calQ_+)$. 
  If $P_0(0,\ldots,0,Y_{\alpha(\calQ_0)},0,\ldots,0)\not=0$, then the unique power 
  series solution has a positive radius of convergence.
}

\bigskip

\Proof We have
$$
  \sum_{\ss A\in\calQ_0\atop\ss \alpha_\ell=\alpha(\calQ_0)} r_A
  = \sum_{-k\leq i\leq k} [Y_i]P(0,Y_{-k},\ldots,Y_k)
  \One\{\, i=\alpha(\calQ_0)\,\} 
$$
and since the nonshifting factors are linear, this 
is 
$$
  P(0,\ldots,0,Y_{\alpha(Q_0)},0,\ldots,0) \, .
$$
The result then follows from Theorem \convergenceDivergence.\hfill\qed

\bigskip

To implement the results of section \fixedref{2.4}, we need an efficient
way to calculate the maximal height of the shifting factors involved
in $q$-algebraic equations. For this, we decompose the shifting part
of the equation as
$$
  \sum_{A\in\calQ_+}r_A A
  = \sum_{-k\leq i\leq k} \sum_{\ss A\in\calQ_+\atop\ss \alpha_\ell=i}
    r_A A \, .
$$
Therefore,
$$
  H(\calQ)
  = \max_{-k\leq i\leq k} \max_{\ss A\in\calQ_+\atop\ss\alpha_\ell=i} i/2a
  = \max_{-k\leq i\leq k} i/\Bigl(2\min_{\ss A\in\calQ_+\atop\ss\alpha_\ell=i} a\Bigr)
  \, .
$$
Since $\sum_{\ss A\in\calQ_+\atop \alpha_\ell=i} r_A A$ corresponds to the part of
$P_+$ which contains $Y_i$ but not $Y_{i+1},\ldots,Y_k$, it corresponds
to
$$
  P_+(z,Y_{-k},\ldots,Y_i,0,\ldots,0)-P_+(z,Y_{-k},\ldots, Y_{i-1},0,\ldots,0)
  \, .
$$
We view this difference as a polynomial $R_i(z)$ in $z$. Then the height
$H(\{\, A\in\calQ_+\,:\,\alpha_\ell=i\,\})$ is $i$ divided by twice the
order of $R(z)$. This leads to the following algorithm.

\bigskip

\hfuzz=2pt
\algorithm{

for $i$ from $-k$ to $k$ by increment of $1$ do

\quad $R_i(z)\leftarrow P_+(z,Y_{-k},\ldots, Y_i,0,\ldots,0)$

\hskip 120pt ${}-P_+(z,Y_{-k},\ldots, Y_{i-1},0,\ldots,0)$

\quad $H_i\leftarrow i/2\ord_z R_i(z)$

$H(\calQ)\leftarrow\max_{-k\leq i\leq k} H_i$

}

\hfuzz=0pt
\bigskip
Then, the crest contains the shifting part of all monomials obtained by 
taking $i$ such that $\ord_z R_i(z)=i/2H(\calQ)$ and taking $[z^{i/2H(\calQ)}]R_i(z)$; it also contains a part obtained from the nonshifting factor $(0;0)$, that
is, the coefficient
$$
  [Y_0]P_0(0,\ldots,0,Y_0,0,\ldots,0)
  = P_0(0,\ldots,0,1,\ldots,0)
$$
where $1$ corresponds to $Y_0=1$. 

We then view
$[z^{i/2H(\calQ)}]R_i(z)$ as a polynomial in $Y_{-k},\ldots,Y_i$. The scopes of
the corresponding $q$-factors are given by the degree of that polynomials
in $Y_i$, and $\ell$ is the total degree of that polynomial. Thus, the
crest polynomial can be calculated as follows.

\bigskip

\algorithm{

$C(z)\leftarrow P_0(0,\ldots,0,Y_0=1,0,\ldots,0)$

for all $i$ such that $\ord_z R_i(z)=i/2H(\calQ)$ do

\quad $T(Y_{-k},\ldots,Y_i)\leftarrow [z^{i/2H(\calQ)}] R_i(z)$

\quad for all the monomials $M(z,Y_1,\ldots,Y_i)$ of $T(Y_{-k},\ldots,Y_i)$

\quad do

\quad\quad $r_A\leftarrow M(1,1,\ldots,1)$

\quad\quad $s(A)\leftarrow \ord_{Y_i} M(Y_{-k},\ldots,Y_i)$

\quad\quad $a\leftarrow i/2H(\calQ)$

\quad\quad $\ell\leftarrow\Deg M(Y_{-k},\ldots,Y_i)$

\quad\quad $C(z)\leftarrow C(z)+r_A s(A) q^{-H(\calQ)a(a-h(\calQ))} z^a t^{\ell-1}$

}

\bigskip

\noindent At the end of this procedure, $C(z)$ is the crest polynomial
$\calC_{q,t}(z)$.

\bigskip

\def\prevs{\the\sectionnumber.\the\snumber }
\def\preveq{(\the\sectionnumber.\the\equanumber)}

\section{Further examples} The purpose of this section is to give
further examples.

\medskip

\noindent{\it Example 1.}
Motivated by combinatorics of lattice paths,
Drake (2009) considers a set of integers $S_0$ which does not contain $0$ and
the equation
$$
  f(z)=1+qzf(z)f(q^2z)+\sum_{j\in S_0}q^{j(j-1)} z^j 
  f(z)f(q^2z)\cdots f(q^{2(j-1)}z) \, .
  \eqno{\equa{DrakeA}}
$$
Though Drake (2009) allows for $S_0$ to be infinite, we restrict $S_0$ 
to be finite
here (extending Theorem \mainTh\ to allow $\calQ$ to be infinite is possible
in some cases, and in a subsequent paper, we will provide such an extension
in the context of the $q$-Lagrange inversion). In our notation, Drake's 
set $S$ is $S_0\cup\{\,0\,\}$.

Equation \DrakeA\ is of the form \qEq\ with $P(z)=1$. Applying $[z^0]$
to both sides of \DrakeA,  we have $f_0=1$. Since $f(z)=(0;0)f(z)$, 
$$
  \calQ=\{\, (0;0), (1;0,2)\,\}\cup \bigl\{\, \bigl(j;0,2,\ldots,2(j-1)\bigr) 
  \,:\, j\in S_0\,\bigr\} \, .
$$
The height of the $q$-factor $(1;0,2)$ is $1$, while that 
of $\bigl(j;0,2,\ldots,2(j-1)\bigr)$ is $(j-1)/j$, which is less than $1$.
Thus, the height of $\calQ$ is $1$; its co-height is $1$; and its crest is
$\widehat\calQ=\{\, (0;0),(1;0,2)\,\}$. This crest has a unique shifting
$q$-factor whose coefficient in equation \DrakeA\ is $r_{(1;0,2)}=q$.
Corollary \mainCor\ yields
$$
  [z^n]f\sim q^{n(n-1)} c q^n = c q^{n^2}
$$
as $n$ tends to infinity. Drake's combinatorial argument, while less general
than ours, is more specific for that equation. Indeed, equation (5) in 
Drake (2009) is our equation \DrakeA, so that his $\tilde r_n^{(S)}(q)$ is
our $f_n$ and his $r_n^{(S)}(1/q)$ is our $q^{-n^2}f_n$. The constant $c$ is
then given in Drake's (2009) Theorem 1 and is
$$
  c=\prod_{i\geq 1} {1\over 1-q^{-2i}} \Bigl(1+\sum_{j\in S_0} q^{-j(2i-1)}\Bigr) 
  \, .
$$

\medskip

\noindent{\it Example 2.} Drake (2009, display (12)) considers also the example
$$
  f(z)=1+zf(z)+qz^2 f(z)f(qz) \, .
  \eqno{\equa{DrakeB}}
$$
He shows that there exist positive constants $\tilde c_0$ and $\tilde c_1$, 
such that
$$
  f_{2m}\sim q^{m^2} \tilde c_0 
  \qquad\hbox{ and }\qquad
  f_{2m+1}\sim q^{m^2+m}\tilde c_1
  \eqno{\equa{DrakeC}}
$$
as $n$ tends to infinity. Note that if $n=2m$, then $m^2=n^2/4$, while if
$n=2m+1$, then $m^2+m=(n^2/4)-1/4$. Thus, setting $c_0=\tilde c_0$ and
$c_1=q^{-1/4}\tilde c_1$, \DrakeC\ can be rewritten as
$$
  f_n\sim q^{n^2/4} c_m
  \eqno{\equa{DrakeD}}
$$
as $n$ tends to infinity, with $m$ being $0$ if $n$ is even and $m$ being $1$
if $n$ is odd.

For the $q$-functional equation \DrakeB, the set of $q$-factors involved is
$$
  \calQ=\{\, (0;0), (1;0), (2;0,1)\,\} \, .
$$
Its height is $1/4$ and co-height $2$. The crest is $\{\, (0;0), (2;0,1)\,\}$
and contains a unique shifting $q$-factor. The corresponding 
coefficient $r_{(2;0,1)}$ in 
equation \DrakeB\ is $q$. Corollary \mainCor\ implies that there exist $c_0$
and $c_1$ such that
$$
  [z^n]f(z)\sim q^{(1/4)n(n-2)} q^{n/2}c_m
   = q^{n^2/4} c_m
$$
where $m$ is the remainder of the division of $n$ by $2$. Therefore, we 
recovered \DrakeD, however with the possibility that $c_0$ or $c_1$ may vanish.
To see that they do not vanish, we
can either use the proof of Corollary \mainCor, or make a direct calculation.
Indeed, since all the coefficients involved in equation \DrakeB\ are
positive, Theorem \mainTh.iii asserts that the function $U_{q,f_0}$ 
in \basicFormula\ has nonnegative coefficients.
The crest polynomial is 
$$
  \calC_q(z)=1-qz^2
$$
and has root $1/\sqrt q$ and $-1/\sqrt q$. We then have, writing $U$ for
$U_{q,f_0}$ in \basicFormula,
$$
  [z^n]{U(z\sqrt q)\over 1-qz^2}
  = {1\over 2} [z^n]\Bigl( {U(z\sqrt q)\over 1-\sqrt q z}
    +{U(z\sqrt q)\over 1+\sqrt q z}\Bigr) \, .
$$
Using singularity analysis, we obtain
$$\eqalign{
  [z^n]{U(z\sqrt q)\over 1-qz^2}
  &{}\sim {1\over 2} \bigl( U(1)q^{n/2}+U(-1)(-1)^n q^{n/2} \bigr) \cr
  &{}\sim q^{n/2} {1\over 2}\bigl( U(1)+(-1)^nU(-1)\bigr) \, . \cr} 
$$
So we have $c_0=\bigl(U_q(1)+U_q(-1)\bigr)/2$ and $c_1=\bigl(U_q(1)-U_q(-1)\bigr)/2$.
Since $U_q$ has nonnegative coefficients, $U_q(1)-U_q(-1)$ does not vanish, and we
have indeed $f_n\sim q^{n^2/4} c_m$ with $c_0$ and $c_1$ not being $0$.

The combinatorial argument in Drake (2009) provides an explicit value for $c_0$
and $c_1$. Indeed, his $\tilde m_n^{(1)}(q)$ is our $f_n$, so that 
his $m_{2n}^{(1)}(1/q)$ is our $q^{-n^2}f_{2n}$ and 
his $m_{2n+1}^{(1)}(1/q)$ is 
our $q^{-(n^2+n)}f_{2n+1}$. Since we have $q^{-n^2}f_{2n}\sim \tilde c_0$ 
and 
$$
  q^{-(n^2+n)}f_{2n+1}
  =q^{1/4} q^{-(2n+1)^2/4}f_{2n+1}
  \sim \tilde c_1 \, ,
$$
as $n$ tends to infinity, it follows from Drake's (2005) Theorem 5 
that $c_0$ is his $\Phi_2(1/q)$ and $c_1$ his $q^{-1/4}\Psi_2(1/q)$, 
that is, with $r=1/q$,
$$
  c_0
  ={1\over (r;r)_\infty (r^2;r^{12})_\infty 
           (r^9;r^{12})_\infty (r^{10};r^{12})_\infty }
$$
and
$$
  c_1
  ={q^{-1/4}\over
     (r;r^2)_\infty^2 (r^4;r^{12})_\infty (r^6;r^{12})_\infty
     (r^{8};r^{12})_\infty (r^{12};r^{12})_\infty } \, .
$$

\medskip

\noindent{\it Example 3.} Gessel (1980) provides also interesting examples
of $q$-functional equations, motivated by the $q$-Lagrange inversion. 
For instance, he considers the equation (Gessel, 1980, equation (10.16))
$$
  f(z)=1+q(1+s)zf(z)+q^3s z^2h(z)h(qz) \, .
  \eqno{\equa{GesselA}}
$$
The $q$-factors involed in this equation are $(0;0)$, $(1;0)$ and $(2;0,1)$. 
The crest is reduced to the $q$-factor $(0;0)$ and the shifting one $(2;0,1)$. 
The corresponding coefficient in equation \GesselA\ is $r_{(2;0,1)}=q^3s$. 
Corollary \mainCor\ yields
$$
  [z^n]f\sim q^{(1/4)n(n-2)} q^{3n/2} s^n c_m
  \sim q^{(n^2/4)+n}s^n c_m
$$
with $m$ being $0$ or $1$ according to the parity of $n$. Again, as in the
previous example, the coefficients $c_0$ and $c_1$ do not vanish when the
coefficients in equation \GesselA\ are positive, that is when $s$ is positive.

\medskip

\noindent{\it Example 4.} Following Garoufalidis (2004), the colored Jones
polynomials for the figure 8 knot are
$$
  J_n(q)=\sum_{0\leq k\leq n-1} q^{nk} 
  \Bigl( {1\over q^{n+1}};{1\over q}\Bigr)_k \Bigl( {1\over q^{n-1}};q\Bigr)_k
$$
with $J_0(q)=1$. Recall the operator $\sigma f(z)=f(qz)$. Garoufalidis (2004)
showed that the generating 
function $J(z)=\sum_{n\geq 0} J_n z^n$ satisfies the $q$-algebraic equation
$$
  C_0 J(z)+z C_1 J(z)+z^2 C_2 J(z)+z^3 C_3 J(z)=0 \, ,
$$
with
$$\eqalign{
  C_0&{}=q\sigma(q^2+\sigma)(q^5-\sigma^2)(1-\sigma) \cr
  C_1&{}=-q^2\sigma^{-1}(1+\sigma)\bigl(q^4-\sigma q^3(2q-1)
       -q^3\sigma^2(q^2-q+1)\cr
     &\hskip 100pt{}+q^4\sigma^3(q-2)+\sigma^4q^4\bigr) (q^3-\sigma^2)
       (1-\sigma) \cr
  C_2&{}=q^7 \sigma^{-1}(1-\sigma)(1+\sigma)(1-q^3\sigma^2)\bigl( q\sigma(q-2)
       +\sigma^2(-1+q-q^2)\cr
     &\hskip 100pt{}-\sigma^3(2q-1)+q\sigma^4\bigr) \cr
  C_3&{}=-q^{10}\sigma(1-\sigma)(1+q^2\sigma)(1-q^5\sigma^2) \, .\cr
  }
$$
The nonshifting $q$-factors come from $C_0$. By considering $C_0$ as a
polynomial in $\sigma$, we see that $\alpha(\calQ_0)=\deg_\sigma C_0=5$. Thus,
we consider $f(z)=\sigma^5J(z)$ as our new function. It solves the
$q$-algebraic equation
$$
  (D_0+zD_1+z^2D_2+z^3D_3)f(z)=0
$$
with $D_i=C_i\sigma^{-5}$. The $D_i$ are Laurent polynomials.
This equation has now $\alpha(\calQ_0)=0$. It 
involves many $q$-factors, all of the form $(a;\alpha_1,\ldots,\alpha_\ell)$
for $a=0,1,2,3$, the $q$-factor $(a;\alpha_1,\ldots,\alpha_\ell)$ coming from
$D_a$. Thus, $\max_{a=i}\alpha_\ell$ is $\deg_\sigma D_i$. For this equation,
the following table
then shows how to calculate efficiently the quantities 
involved in theorem \mainTh.

\bigskip

\setbox1=\vbox{\halign{#\hfill\quad&\hfill$#$\quad&\hfill$#$\quad&\hfill$#$\quad&\hfill$#$\cr
$i={}$ & 0 & 1 & 2 & 3 \cr
leading term of $C_i$ & q\sigma^5 & -q^6\sigma^7 &q^{11}\sigma^7 &-q^{17}\sigma^5 \cr

leading term of $D_i$ & q         & -q^6 \sigma^2& q^{11}\sigma^2 & -q^{17}\cr
corresponding $r_A$A&   q(0;0)    & -q^6(1;2)    &q^{11}(2;2)   &-q^{17}(3;0)\cr
$H(A)$              &\hbox{--}    & 1            & 1/2          & 0\cr
$r_A$               & q           & -q^6         & q^{11}       &-q^{17}\cr}}

\centerline{\box1}
\bigskip

From this table, we read that the height of the set of $q$-factors involved
in the equation for $f$ is $1$ and the co-height is $1$. The crest is
$$
  \widehat\calQ=\{\, (0;0), (1;2)\,\} \, .
$$
Corollary \mainCor\ yields
$$
  [z^n]f(z)\sim c q^{n(n-1)} q^{5n}
$$
as $n$ tends to infinity. Since $f_n=q^{5n}J_n$, we conclude that
$$
  J_n\sim c q^{n(n-1)}
  \eqno{\equa{JonesA}}
$$
as $n$ tends to infinity.

Since the assumptions for the last assertion of Theorem \mainTh\ does not hold
in this example, it is possible that $c$ is $0$. This example is instructive
since the asymptotic behavior of $J_n$ can be obtained directly from its
expression. Indeed, we have
$$\eqalign{
  J_n
  &{}=\sum_{0\leq k\leq n-1} q^{nk} 
    {1-1/q^n\over (1-1/q^{n-k})\cdots (1-1/q^{n+k})} \cr
  &{}=\sum_{0\leq j\leq n-1} q^{n(n-1-j)} 
    {1-1/q^n \over (1-1/q^{j+1})\cdots (1-1/q^{2n-1-j})} \, .\cr
  }
$$
Therefore, isolating the term for which $j=0$,
$$\displaylines{\qquad
  J_n=q^{n(n-1)} {1-1/q^n\over (1/q;1/q)_{2n-1}}
  \hfill\cr\hfill
  {}+ q^{n(n-1)} (1-1/q^n) \sum_{1\leq j\leq n-1} 
  {q^{-nj}\over (1/q^{j+1};1/q)_{2(n-j)-1}} \, .
  \qquad\cr}
$$
Since $|q|>1$, we have $0\leq 1-1/|q^j|\leq |1-1/q^j|$ for any nonnegative
integer $j$, and therefore,
$$\eqalign{
  \Bigl|\sum_{1\leq j\leq n} {q^{-nj}\over (1/q^{j+1};1/q)_{2(n-j)-1}}\Bigr|
  &{}\leq \sum_{1\leq j\leq n} {|q|^{-nj} \over (1/|q|;1/|q|)_\infty} \cr
  &{}\leq {1\over (1/|q|;1/|q|)_\infty} {|q|^{-n}\over 1-|q|^{-n}} \, . \cr
  }
$$
Consequently,
$$
  J_n\sim {q^{n(n-1)}\over (1/q;1/q)_\infty }
$$
as $n$ tends to infinity. We see that \JonesA\ is in fact sharp.

\bigskip

\noindent{\it Example 5.} To illustrate their results, Cano and Fortuny Ayuso
(2012) consider the $q$-algebraic equation given by the polynomial
$$
  4Y_1^4-9Y_0^2Y_1Y_2+2Y_0^3Y_2 +{z\over q^4}Y_0Y_2
  -z^3Y_0^4Y_5^2-{z^3\over q^4} Y_2 -z^3Y_0+z^5 \, .
  \eqno{\equa{CFAa}}
$$
For this example, we need the imaginary unit $\imath$ such that $\imath^2=-1$.
Considering the particular value $q=4$, they show among other things that 
this $q$-algebraic equation has a unique 
solution $f(z)=z^2+42\sqrt{2} \imath z^{7/2}+o(z^{7/2})$
and that $f_n$ is of $4$-Gevrey order $3/34$. The following result illustrates
how Cano and Fortuny Ayuso's results can be fruitfully used with ours to obtain
information on the asymptotic behavior of the coefficients of generalized 
diverging power series of nearly any $q$-algebraic equation.

\Proposition{\label{CFA}
  For all but countably many values of $q$ of modulus greater than $1$, 
  \CFAa\ has a solution a power series $f(z)$ in $z^{1/2}$ such that
  $$
    f(z)=z^2+\rho z^{7/2}+o(z^{7/2})
  $$
  with $\rho^2=-q(q^2-2)(2q-1)(2q+1)$ and
  $$
    f_n\sim c_{q,m}\, q^{(3n^2-63n)/68}
    \Bigl( {-2\over \rho}\Bigr)^{n/17} \, .
  $$
  for $m=n \ \mod 17$ and $0\leq m<17$.
}

\bigskip

In particular, Proposition \CFA\ implies that $f_n$ has $4$-Gevrey order
$3/34$, as indicated by Cano and Fortuny Ayuso (2012).

\bigskip

Our proof shows that the conclusion of Proposition \CFA\ is valid for the
specific value $q=4$. When $q=4$ we have $\rho^2=-42^2\times 2$. So, we
may take $\rho=42 \imath\, \sqrt 2$, which leads to 
$$
  f_n\sim c_m \imath^{n/17} \, {2^{(3n^2-64n)/34}\over 21^{n/17}}
  \eqno{\equa{CFAaa}}
$$
as $n$ tends to infinity. Since $\imath^{(n+17)/17} =\imath\,\imath^{n/17}$,
we may rewrite \CFAaa\ as
$$
  f_n\sim \tilde c_m {2^{(3n^2-64n)/34}\over 21^{n/17}}
$$
with now $m=n\ \mod\, 68$ and the additional constraint that $\tilde c_{m+17}
=\imath\, \tilde c_m$ for $0\leq m<51$. Using {\tt Maxima}, Jos\'e Cano 
computed exactly the first 390 coefficients of the solution and kindly 
shared and allowed us to use the
result of this computation. An inspection of the coefficients 
reveals that $f_n$ is real if $n$ is even and purely imaginary if $n$ is odd, 
and that the sign of $f_n$ is the opposit of that of $f_{n+34}$, 
confirming the relation $\tilde c_{m+17}=\imath\,\tilde c_m$. The following 
plot shows $\log |f_n {21^{n/17}/ 2^{(3n^2-64n)/34}}|$.

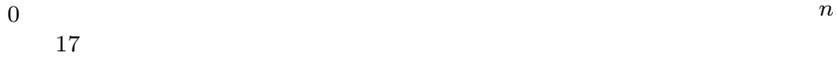
\setbox1=\hbox to 300pt{\hskip 16pt\vbox to 160pt{
  \eightpoints
  \hsize=300pt
  \kern 160pt
  \includegraphics{qFunctEq_fig.eps}
  \anotePT{-9}{94}{$0$}
  \anotePT{270}{96}{$n$}
  \anotePT{7.4}{84}{$17$}
  \anotePT{-19}{0}{$-45$}
  \anotePT{-13}{141}{$25$}
  \vfill
  \vss
}\hfill}    

\box1

\bigskip

\noindent
This plots confirms our result, and could be used to evaluate the 
coefficients $c_m$ numerically. After Theorem \mainTh, we 
made the remark that we can provide an estimate for $\Theta$. The positivity
of $\Theta$ implies that the convergence of
$$
  \Bigl({-\rho\over 2}\Bigr)^{n/17} q^{-(3n^2-63n)/68} f_n-\tilde c_m
$$
to $0$ is exponentially fast. This is confirmed by the plot which exhibits
a complete stability for $n$ greater than about $70$. While this plot 
supports that some $\tilde c_m$ do not vanish, it is unclear if those of 
magnitude less than $10^{-10}$ say are zero or not.

\bigskip

\Proof Cano and Fortuny-Ayuso's (2012) method provides the crucial 
information that
for $q=4$, the equation has a solution which is a power series in $z^{1/2}$. 
Equipped with this information, we try in general to find a solution as a
power series in $z^{1/2}$, which leads to the change of function $g(z)=f(z^2)$
and to define $r$ such that $r^2=q$; thus $r$ is defined up to a sign.
Substituting $z^2$ for $z$ in the equation given by \CFAa, we see that if $f$
is a solution, then $g$ solves the equation given by
$$
  P=4Y_1^4-9Y_0^2Y_1Y_2+2Y_0^3Y_2+{z^2\over r^8} Y_0Y_2 - z^6 Y_0^4Y_5^2
  -{z^6\over r^8} Y_2 -z^6Y_0+z^{10} \, .
$$
The study of this equation is made with the algorithms described in the 
previous section, and the {\tt Maple} code that we used is in the appendix
to this paper; this code is a very simple and primitive implementation,
easy to read and to adapt to other equations; this is not a polished package.

The equation given by $P$ is not reduced since the nonshifting 
part contains the nonlinear
terms $Y_1^4$, $Y_0^2Y_1Y_2$ and $Y_0^3Y_2$. We apply the reduction algorithm. 
In fact we need to apply the reduction step 4 times to obtain that
$g_0=g_1=g_2=g_3=0$ and $g_4=1$. Setting $g(z)=z^4+z^5h(z)$, the function $h$
solves an equation given by a polynomial made of 41 monomials. This equation
is still not reduced. We then apply the reduction procedure 3 more times,
to obtain $g_5=g_6=0$ and $g_7$ solves
$$
  g_7^2=r^2(2-r^4)(2r^2-1)(2r^2+1) \, .
$$
Thus, $g_7$ is to be chosen up to a sign, and we may set $g_7=\rho$. 
Setting $g(z)=z^2+\rho z^7+z^8h(z)$, we
obtain that $h_0=0$. Thus, to obtain an estimate as accurate as possible, 
following the remark after Corollary \qGevrey, we apply the reduction step 3 
more times, to obtain that $g_8=g_9=0$ and that $g_{10}$ solves
$$
  (r^6+1)g_{10} +r^2(6-18r^4+2r^6-9r^7-9r^{10}+16r^{11}) =0\, .
$$
Since $|r|>1$, this equation has a nonvanishing solution
$g_{10}$ except when $r$ is one of the $11$ roots of the polynomial 
$6-18r^4+2r^6-9r^7-9r^{10}+16r^{11}$. So we set
$$
  g(z)=z^4+\rho z^7+z^{10}h(z)
$$
and obtain that $h_0\not=0$ and that $h$ satisfies a $r$-algebraic equation.
After simplification, the polynomial corresponding to that equation is made
of $397$ monomials. The nonshifting factors are given by the polynomial
$$
  P_0=\rho (r^6Y_2+Y_0) \, .
$$
Thus, writing $\calQ$ for the set of $r$-factors involved in this equation
for $h$, we have $\alpha(\calQ_0)=2$. The shifting $r$-factors are contained
in a polynomial $P_+$ which is made of $292$ monomials. An application of the
algorithm given in the previous section shows that $\alpha(\calQ_+)=5$. To 
check the existence condition \algoExistenceConditionA, a computation with
{\tt Maple} shows that
$P(0,Y_0,\ldots,Y_5)$ is proportional to 
$$
  Y_0+r^6Y_2+r^2(6-18r^4+2r^6-9r^7-9r^{10}+16r^{11}) \, .
$$
Thus, this existence condition is that
$$
  r^{6+2n}+1+6r^2-18r^6+2r^8-9r^9-9r^{12}+16r^{13}\not=0
$$
for all $n$. For a given $n$, this condition fails for at most $(6+2n)\vee 13$
values of $r$. Therefore, this condition fails for at most countably many $r$. 
In particular, if $r=2$, this condition becomes
$$
  64\cdot 4^n+88\,985\not=0
$$
which is of course verified.

Having verified the condition, we need to transform the equation to one for
which $\alpha(\calQ_0)$ vanishes so that we can apply Theorem \mainTh. This
is done by setting $h(z)=k(z/r^2)$. For the new equation in $k$, the algorithm
of the previous section yields that the height of the equation is
$3/34$ and the co-height is $17$. The crest is given by the following part
of the equation,
$$
  \rho r^6 Y_0-2r^{64}z^{17}Y_3 \, .
$$
In particular, the crest has a unique shifting $r$-factor, $-2r^{64}z^{17}Y_3$.
We then apply Corollary \mainCor\ to obtain
$$
  k_n
  \sim c_{r,m} r^{3n(n-17)/34} \Bigl( -{2r^{64}\over \rho r^6}\Bigr)^{n/17}
  \sim c_{r,m} r^{(3n^2+65n)/34} \Bigl(-{2\over \rho}\Bigr)^{n/17}
$$
as $n$ tends to infinity. Given our change of function, 
we have for any $n$ at least $10$,
$$
  f_n=g_n=h_{n-10} = k_{n-10}/r^{2(n-10)}
$$
The result follows after some calculations.

\bigskip

\def\prevs{\the\sectionnumber.\the\snumber }
\def\preveq{(\the\sectionnumber.\the\equanumber)}

\section{Proof of the Theorems}
The main tool is the recursion
that the coefficients of the solution must obey. We consider a $q$-algebraic
equation \qEq. We write $P(z)=\sum_{0\leq i\leq p}P_i z^i$. 
Applying $[z^n]$ to both sides of the equation, we obtain for any $n\geq 0$,
$$
  \sum_{A\in\calQ} r_A \sum_{n_1+\cdots+n_\ell=n-a} 
  q^{\alpha_1n_1+\cdots+\alpha_\ell n_\ell} f_{n_1}\cdots f_{n_\ell} 
  +P_n =0 \, .
  \eqno{\equa{fnRec}}
$$

\def\prevs{\the\sectionnumber.\the\subsectionnumber.\the\snumber }
\def\preveq{(\the\sectionnumber.\the\subsectionnumber.\the\equanumber)}

\subsection{Proof of Theorem \existence}
In \fnRec\ we distinguish according
to the shifting and nonshifting $q$-factors. Since the nonshifting ones
are linear we obtain
$$\displaylines{\quad
  -\sum_{A\in\calQ_0}r_A q^{\alpha_1 n} f_n
  \hfill\cr\hfill
  = \sum_{A\in\calQ_+} r_A \sum_{n_1+\cdots+n_\ell=n-a}
  q^{\alpha_1n_1+\cdots+\alpha_\ell n_\ell}
  f_{n_1}\cdots f_{n_\ell}
  {}+P_n\, .\quad\equa{fnRecA}\cr}
$$
By definition all $q$-factors in $\calQ_+$ are shifting and therefore have
$a\geq 1$. Thus, we see that under \existenceConditionA, we can calculate $f_n$
inductively once $f_0$ is determined.

\subsection{Proof of Theorem \convergenceDivergence}
Let $L=\max_{A\in\calQ} \ell$. We first consider the sequence defined 
inductively by $g_0=1$ and for any $n\geq 1$,
$$
  g_n=\sum_{n_1+\cdots+n_L=n-1} g_{n_1}\cdots g_{n_L} \, .
  \eqno{\equa{gnBasic}}
$$
This sequence is nonnegative. Furthermore, considering the tuple 
$(n_1,\ldots,n_L)=(n-1,0,\ldots,0)$ in \gnBasic, we see 
that $g_n\geq g_{n-1}$.

\Lemma{\label{gnGFCv}
  The generating function $g(z)=\sum_{n\geq 0} g_n z^n$ has a positive
  and finite radius of convergence.
}

\bigskip

\Proof The proof of this lemma is inspired by an argument in F\"urlinger 
and Hofbauer (1985). Given recursion \gnBasic, considering $g(z)$ as a formal
power series, we have
$$
  g(z)=1+zg(z)^L \, .
$$
Thus,
$$
  g(z^{L-1})=1+z^{L-1} g(z^{L-1})^L \, .
  \eqno{\equa{gnBasicA}}
$$
Set $k(z)=zg(z^{L-1})$. Then $k(z)$ is a formal power series and
$k(z)=\sum_{n\geq 1} k_n z^n=z\sum_{n\geq 0} g_n z^{n(L-1)}$. In particular
the sequence $(k_n)$ is nonnegative. Multiplying both sides of \gnBasicA\ by
$z$, we have
$$
  k(z)=z+k(z)^L \, .
  \eqno{\equa{gnBasicB}}
$$
Setting $y=k(z)$, \gnBasicB\ becomes $z=y(1-y^{L-1})$. The 
relation $k(z)=\sum_{n\geq 1} k_n z^n$ then becomes the relation between 
power series,
$$
  y=\sum_{n\geq 0} k_n \bigl( y(1-y^{L-1})\bigr)^n \, .
$$
Since $(k_n)$ is a nonnegative sequence, since
$0\leq 1-y^{L-1}\leq 1$ whenever $0\leq y\leq 1$, this relation leads to
$$
  y\geq \sum_{n\geq 0} k_n y^n =k(y)\geq 0\, ,
$$
in the range $0\leq y\leq 1$. This shows that the radius of convergence 
of $k(z)$ is positive. Therefore,
the radius of convergence of $g$ is also positive.\hfill\qed

\bigskip

Since we assume that $\alpha(\calQ_0)\geq\alpha(\calQ_+)$, the discussion at
the end of subsection \fixedref{2.2} shows that 
we can assume without any loss of generality that $\alpha(\calQ_0)=0$ so
that $\alpha(\calQ_+)\leq 0$.

\Lemma{\label{simpleBound}
  Consider a reduced $q$-algebraic equation whose $q$-factors have been
  collected. If $\alpha(\calQ_0)=0$ and \existenceConditionA\ holds, then
  $$
    \inf_{n\in\NN}  \Bigl|\sum_{A\in\calQ_0} r_A q^{\alpha_1 n}\Bigr| \not= 0 \, .
  $$
}

\Proof Since \existenceConditionA\ holds, $\sum_{A\in\calQ_0}r_Aq^{\alpha_1n}$ does
not vanish for any $n$. Since $\alpha(\calQ_0)=0$, all the $\alpha_\ell$
pertaining to a nonshifting $q$-factor are nonpositive. Thus,
$$
  \lim_{n\to\infty} \Bigl|\sum_{A\in\calQ_0}r_A q^{\alpha_1 n}\Bigr|
  =\sum_{\ss A\in\calQ_0\atop\ss \alpha_\ell=0} r_A \, .
  \eqno{\equa{simpleBoundA}}
$$
Since the $q$-algebraic equation is reduced and its $q$-factors have been
collected, there is a unique nonshifting $q$-factor with $\alpha_\ell=0$,
and $r_{(0;0)}\not=0$. Hence the right hand side of \simpleBoundA\ 
is $r_{(0;0)}$ and does not vanish.\hfill\qed

\bigskip

Let $(g_n)$ be defined in \gnBasic, with $g_0=1$. To prove 
Theorem \convergenceDivergence, it suffices to prove that there 
exists some positive $c$ and some $\gamma$ greater than $1$ such 
that $|f_n|\leq c\gamma^n g_n$. The proof is by induction. We 
set $c=1\vee\max_{0\leq n\leq p} |f_n|$. We will see how to choose $\gamma$ 
afterwards. Given Lemma \simpleBound,
$$
  r=\inf_{n\in\NN} \Bigl| \sum_{A\in\calQ_0} r_A q^{\alpha_1 n}\Bigr|
$$
is positive. 

Let $n$ be some integer greater than $p$. Suppose that 
we proved that $|f_i|\leq c\gamma^i g_i$ for 
any $0\leq i\leq n-1$. Consider the recursion \fnRecA. If $A$ is 
in $\calQ_+$ then all the $\alpha_i$ are 
nonpositive and $|q^{\alpha_1n_1+\cdots+\alpha_\ell n_\ell}|\leq 1$ since $|q|>1$.
Thus, \fnRecA\ yields
$$
  r |f_n|
  \leq  \sum_{A\in\calQ_+} |r_A| \sum_{n_1+\cdots+n_\ell=n-a}
    |f_{n_1}|\cdots |f_{n_\ell}| \, .
  \eqno{\equa{fnBoundA}}
$$
Using the induction hypothesis and the notation $L=\max_{A\in\calQ}\ell$, 
this is at most
$$
  c^L\gamma^{n-1} \sum_{A\in\calQ_+} |r_A|\sum_{n_1+\cdots+n_\ell=n-a}
  g_{n_1}\cdots g_{n_\ell} \, .
$$
Since $(g_n)$ is nondecreasing, we see that whenever $a$ is positive, 
$g_{n_1}\leq g_{n_1+a-1}$. Since a tuple $(n_1,\ldots,n_\ell)$
may be seen as tuples $(n_1,\ldots,n_\ell,0,\ldots,0)$ of length $L$ and
$g_0=1$,
$$
  \sum_{n_1+\cdots+n_\ell=n-1} g_{n_1}\cdots g_{n_\ell}
  \leq \sum_{n_1+\cdots+n_L=n-1} g_{n_1}\cdots g_{n_L} \, .
  \eqno{\equa{fnBoundB}}
$$
Thus, \fnBoundA\ yields
$$
  r |f_n|
  \leq c^L \gamma^{n-1}\sum_{A\in\calQ_+} |r_A| g_n \, .
$$
We choose $\gamma\geq c^{L-1}\sum_{A\in\calQ_+}|r_A|/r$ to 
obtain $|f_n|\leq c\gamma^n g_n$ for all $n$ greater than $p$. If $\gamma$ is
large enough, this inequality holds for all integers $n$.

\subsection{Proof of Theorem \mainTh}
The sequence
$$
  d_n=H(\calQ)n(n-h(\calQ))
$$
is key to our proof. For simplicity of notation, we will write $H$ 
for $H(\calQ)$ and $h$ for $h(\calQ)$, so that $d_n=Hn^2-Hhn$. We also agree
that $d_n=0$ if $n$ is negative. Recall that a sequence like $(d_n)$ is 
strictly convex if $(d_{n+1}-d_n)$ is an increasing sequence.

Throughout this proof, $(n_1,\ldots,n_\ell)$ denotes a tuple of nonnegative
integers.

\Lemma{\label{dnConvex}
  The sequence $(d_n)$ is strictly convex.
}

\bigskip

\Proof Clearly, $d_{n+1}-d_n=(2n+1)H-Hh$ is increasing in $n$.\hfill\qed

\bigskip

For a $q$-factor $A=(a;\alpha_1,\ldots,\alpha_\ell)$, set
$$
  D_A(n_1,\ldots,n_\ell)
  =d_{n_1+\cdots+n_\ell+a}-(d_{n_1}+\cdots+d_{n_\ell}
  +\alpha_1n_1+\cdots+\alpha_\ell n_\ell) \, .
$$

The following definition is stated to set the notation and recall some
known terminology (see MacDonald, 1995, \S I.1 for raising operators).

\Definition{\label{operators}
  Given two positive integers $i<j$, the following operators act on 
  tuples of length at least $j$:

  \medskip

  \noindent (i) the transposition $\tau_{i,j}$ permutes the $i$-th 
  and $j$-th entries.

  \medskip

  \noindent (ii) the raising operator $R_{i,j}$ decreases the $i$-th 
  entry by $1$ and increasing the $j$-th entry by $1$.
}

\bigskip

For instance $R_{3,5}(1,2,3,4,5,6)=(1,2,2,3,4,6,6)$. Note that transpositions
and raising operators leave invariant the sum of the entries of a tuple.

\Lemma{\label{DTransform}
  Let $A=(a;\alpha_1,\ldots,\alpha_\ell)$ be a $q$-factor and 
  let $1\leq i<j\leq\ell$ be some integers. Then
  $$
    D_A\circ\tau_{i,j}(n_1,\ldots,n_\ell)=D_A(n_1,\ldots,n_\ell)
  +(\alpha_j-\alpha_i)(n_j-n_i)
  $$
  and,
  $$\displaylines{\qquad
    D_A\circ R_{i,j} (n_1,\ldots,n_\ell)=D_A(n_1,\ldots,n_\ell)
    +d_{n_i}-d_{n_i-1}
    \hfill\cr\hfill
    {}-(d_{n_j+1}-d_{n_j})+\alpha_i-\alpha_j \, .\qquad\cr}
  $$
}

\Proof Both assertions follow from some elementary calculation.\hfill\qed

\bigskip

Because it is a cumbersone quantity in our coming calculations, 
for a $q$-factor $A=(a;\alpha_1,\ldots,\alpha_\ell)$ and a positive 
integer $k$ at most $\ell$ we set 
$$
  \theta(A,1)=Ha^2+(Hh-\alpha_\ell)a \, .
$$
and for $2\leq k\leq \ell$,
$$\displaylines{\quad
  \theta(A,k)=H(a+k-1)^2+(Hh-\alpha_\ell)(a+k-1)
  \hfill\cr\hfill
  {}+(k-1)H(1-h)
  +(\alpha_{\ell-k+1}+\cdots+\alpha_{\ell-1}) \, .
  \quad\cr}
$$

\Lemma{\label{Dmin}
  Let $A$ be a $q$-factor. Let $n$ be a nonnegative integer with $n\geq a$.
  The smallest value of $D_A(n_1,\ldots,n_\ell)$ when $n_1+\cdots+n_\ell=n-a$
  and $\sharp\{\, i\,: n_i>0\,\}=k$ is
  $$
    n\bigl(2Ha-\alpha_\ell+2H(k-1)\bigr)
    -\theta(A,k) \, ,
  $$
  and is achieved at the unique tuple 
  $(0,\ldots,0,\underbrace{1,\ldots,1}_{k-1},n-a-k+1)$.
}

\bigskip

\Proof Consider a tuple $(n_1,\ldots,n_\ell)$ whose entries sum to $n-a$.
If $i<j$ and $n_i>n_j$, we apply a transposition $\tau_{i,j}$ so that 
Lemma \DTransform\ yields
$$
  D_A\circ\tau_{i,j}(n_1,\ldots,n_\ell)\leq D_A(n_1,\ldots,n_\ell) \, .
$$
Thus, if we are given $(n_1,\ldots,n_\ell)$ up to a permutation, the smallest
value of $D_A(n_1,\ldots,n_\ell)$ is achieved when $n_1\leq\ldots\leq n_\ell$. 

Consider such an ordered tuple. If $n_i>1$ and $i<j\leq\ell$, Lemma
\dnConvex\ yields
$$
  d_{n_i}-d_{n_i-1}-(d_{n_j+1}-d_{n_j}) \leq 0 \, ,
$$
and, since $\alpha_i\leq\alpha_j$, we also have $\alpha_i-\alpha_j\leq 0$. 
Therefore, Lemma \DTransform\ yields
$$
  D_A\circ R_{i,j}(n_1,\ldots,n_\ell) \leq D_A(n_1,\ldots,n_\ell) \, .
$$

To prove the lemma, we start with an arbitrary tuple $(n_1,\ldots,n_\ell)$ 
with exactly $k$ positive entries. We order it, and
apply several raising operators $R_{i,j}$, $i<j$, to bring it to the
form $(0,\ldots,0,1,\ldots,1,n-a-k+1)$. Each time, $D_A$ decreases.

It follows that the minimum of $D_A$ over the given tuples is achieved at the
unique tuple $(0,\ldots,0,1,\ldots 1,n-a-k+1)$ and is
$$
  d_n-(k-1)d_1 -d_{n-a-k+1}-\bigl(\alpha_{\ell-k+1}+\cdots+\alpha_{\ell-1}
  +(n-a-k+1)\alpha_\ell\bigr) \, .
$$
The conclusion of the lemma then follows from the definition of $d_n$ and an
elementary calculation.\hfill\qed

\bigskip

Lemma \Dmin\ has the following consequences.

\Lemma{\label{DminCor}
  There exist some positive $\Theta$ and $\theta^*$
  such that

  \medskip
  \noindent (i) for any $A$ in $\widehat\calQ$,
  $$
    \min_{\ss n_1+\cdots+n_\ell=n-a\atop\ss\sharp\{\, i\,:\, n_i>0\,\}\geq 2}
    D_A(n_1,\ldots,n_\ell) \geq \Theta n-\theta^* \, ;
    \eqno{\equa{DminCorA}}
  $$
  
  \medskip
  \noindent (ii) for any $A$ in $\calQ\setminus\widehat\calQ$,
  $$
    \min_{n_1+\cdots+n_\ell=n-a}D_A(n_1,\ldots,n_\ell)\geq \Theta n-\theta^*
    \, ;
    \eqno{\equa{DminCorB}}
  $$

  \medskip
  \noindent (iii) If only the $i$-th entry of $(n_1,\ldots,n_\ell)$ does
  not vanish,
  $$
    D_A(0,\ldots,0,n-a,0,\ldots,0)
    =n(2aH-\alpha_i)-Ha^2-Hha+\alpha_i a\, .
  $$

}

\Proof Set $\theta^*=\max_{A\in\calQ}\max_{1\leq k\leq \ell} \theta(A,k)$
where the $\ell$ in the inner maximum pertains to the $A$ in the outer one.
We take $\Theta$ to be
the smallest of $2H$, $\min_{A\in\calQ\setminus\widehat\calQ}(2Ha-\alpha_\ell)$ and
$1$.

\noindent (i) If $A$ is in $\widehat\calQ$, then $2Ha-\alpha_\ell=0$. 
Thus, Lemma \Dmin\ 
implies that for $A$ in $\widehat\calQ$ and 
if $(n_1,\ldots,n_\ell)$ has at least $2$ positive entries, then
$$
  D_A(n_1,\ldots,n_\ell)\geq 2nH -\theta^* \,.
$$

\noindent (ii)
If $A$ is in $\calQ\setminus\widehat\calQ$, then $2Ha-\alpha_\ell$ is positive.
A tuple $(n_1,\ldots,n_\ell)$ with $n_1+\cdots+n_\ell=n-a$ has at least one
positive entry. Lemma \Dmin\ implies that $D_A(n_1,\ldots,n_\ell)$ is
at least $n(2Ha-\alpha_\ell)-\theta^*$. 

\noindent (iii) This is $d_n-d_{n-a}-\alpha_i(n-a)$.\hfill\qed

\bigskip

Recall that throughout the paper we assume $|q|>1$. We define
$$
  g_n=q^{-d_n} f_n \,.
  \eqno{\equa{gnDef}}
$$
It is convenient to agree that $g_n=0$ if $n$ is negative. We write
$$
  P(z)=\sum_{0\leq i\leq p} P_i z^i
$$
the polynomial involved in \qEq.

\Lemma{\label{gnRecurence}
  For any nonnegative integer $n$,
  $$\displaylines{\qquad
    -\sum_{A\in\calQ_0} r_A q^{\alpha_1n} g_n
    = q^{-d_n} P_n
    \hfill\cr\hfill
    {}+\sum_{A\in\calQ_+} r_A \sum_{n_1+\cdots+n_\ell=n-a}
    q^{-D_A(n_1,\ldots,n_\ell)} g_{n_1}\cdots g_{n_\ell} \, . 
    \qquad\equa{gnRecA}\cr}
  $$
}

\Proof
Multiplying both sides of \fnRecA\ by $q^{-d_n}$, 
the result follows after some simple calculations.\hfill\qed

\bigskip

\Lemma{\label{gnGeom}
  There exist some positive $c$ and $G$ such that $|g_n|\leq c\, G^n$ for any
  $n\geq 0$.
}

\bigskip

\Proof The proof is by induction. We write 
$D_A(0,\ldots,n-a,\ldots 0)_i$ for $D_A(0,\ldots,0,n-a,0,\ldots,0)$ 
where $n-a$ is in the $i$-th component.
We set $c=1\vee \max_{0\leq i\leq p} |g_i|$ and we will see how to 
determine $G$. Provided that $G\geq 1$, we have $|g_i|\leq c\,G^i$ for any 
$0\leq i\leq p$. Let $n$ be an integer greater than $p$. Assume that we 
have $|g_i|\leq c\, G^i$ for any $i\leq n-1$.
We split the second term on the right hand side of \gnRecA\ as the sum of
$$\eqalign{
  V_{1,n}&{}=\sum_{A\in\calQ_+\cap\widehat\calQ} r_A 
          \sum_{\ss 1\leq i\leq\ell \atop\ss\alpha_i=\alpha_\ell}
          q^{-D_A(0,\ldots,n-a,\ldots,0)_i} g_{n-a} g_0^{\ell-1} \, ,\cr
  V_{2,n}&{}=\sum_{A\in\calQ_+\cap\widehat\calQ} r_A 
          \sum_{\ss 1\leq i\leq\ell\atop\ss \alpha_i\not=\alpha_\ell}
          q^{-D_A(0,\ldots,n-a,\ldots 0)_i} g_{n-a}g_0^{\ell-1} \, , \cr
  V_{3,n}&{}=\sum_{A\in\calQ_+\cap\widehat\calQ} r_A 
          \sum_{\ss n_1+\cdots+n_\ell=n-a\atop\ss \sharp\{ i:n_i>0\}\geq 2}
          q^{-D_A(n_1,\ldots,n_\ell)} g_{n_1}\cdots g_{n_\ell} \, , \cr
  V_{4,n}&{}=\sum_{A\in\calQ_+\setminus\widehat\calQ} r_A
          \sum_{n_1+\cdots+n_\ell=n-a} q^{-D_A(n_1,\ldots,n_\ell)}
          g_{n_1}\cdots g_{n_\ell} \, . \cr}
$$

Consider $V_{1,n}$. Let $A$ be in $\calQ_+\cap\widehat\calQ$. 
If $\alpha_i=\alpha_\ell$, then Lemma \DminCor.iii yields
$$
  D_A(0,\ldots,n-a,\ldots,0)_i = Ha(a-h) \, .
$$
Thus,
$$
  V_{1,n}
  =\sum_{A\in\calQ_+\cap\widehat\calQ} r_A s(A) q^{-Ha(a-h)} g_{n-a} g_0^{\ell-1}
  \, . \eqno{\equa{gnGeomA}}
$$
Using the induction hypothesis, if $A\in\calQ_+\cap\widehat\calQ$, 
then $a\geq 1$ and $|g_{n-a}|\leq c\, G^{n-1}$. Thus, introducing
$$
  C_1=\sum_{A\in\calQ_+\cap\widehat\calQ} |r_A| s(A)|g_0|^{\ell-1} \, ,
$$
we have
$$
  |V_{1,n}|\leq C_1 c\, G^{n-1} \, .
$$

We now consider $V_{2,n}$. If $A$ is $\calQ_+\cap\widehat\calQ$ and
$\alpha_i\not=\alpha_\ell$, then Lemma \DminCor.iii yields
$$
  D_A(0,\ldots,n-a,\ldots,0)_i
  = n(\alpha_\ell-\alpha_i)+Ha(a-h)+(\alpha_i-\alpha_\ell)a \, .
$$
Moreover, if $\alpha_i\not=\alpha_\ell$, then $\alpha_\ell-\alpha_i\geq 1$.
Hence,
$$
  D_A(0,\ldots,n-a,\ldots,0)_i
  \geq n+Ha(a-h)+(\alpha_i-\alpha_\ell)a \, .
  \eqno{\equa{gnGeomB}}
$$
Therefore, using the induction hypothesis, there exists $C_2$ such that
$$
  |V_{2,n}|\leq |q|^{-n} C_2 c\, G^{n-1} \, .
$$

We consider now $V_{3,n}$. In what follows, we will write $\ell^*$ for 
$\max\{\, \ell\,:\, A\in\calQ\,\}$. If $A$ is in $\calQ_+\cap\widehat\calQ$ 
and two $n_i$ are positive, Lemma \DminCor.i yields
$$
  D_A(n_1,\ldots,n_\ell)\geq \Theta n-\theta^* \, .
$$
Moreover, if $n_1+\cdots+n_\ell=n-a$, then the largest $n_i$ is at most
$n-a$, which is at most $n-1$ since $A$ is shifting. Therefore, using the 
induction hypothesis, and bounding by $n^\ell$ the number of 
tuples $(n_1,\ldots,n_\ell)$ with at least two postive entries 
and $n_1+\cdots+n_\ell=n-a$, we obtain
$$\eqalign{
  |V_{3,n}|
  &{}\leq \sum_{A\in\calQ_+\cap\widehat\calQ} |r_A| 
   \sum_{\ss n_1+\cdots+n_\ell=n-a\atop\ss \sharp\{i:n_i>0\}\geq 2}
   |q|^{-n\Theta+\theta^*} c^\ell G^{n-a} \cr
  &{}\leq\sum_{A\in\calQ_+\cap\widehat\calQ} |r_A| |q|^{-\Theta n+\theta^*}
    c^\ell G^{n-1} n^\ell \cr
  &{}\leq |q|^{-\Theta n} n^{\ell^*} \sum_{A\in\calQ_+\cap\widehat\calQ} |r_A|
    |q|^{\theta^*} c^\ell G^{n-1} \, . \cr}
$$
Since the sequence $(|q|^{-\Theta n}n^{\ell^*})$ is bounded, there exists $C_3$
such that
$$
  |V_{3,n}|\leq C_3 c^{\ell^*} G^{n-1} \, .
$$

Finally, if $A$ is in $\calQ_+\setminus\widehat\calQ$, then Lemma \DminCor.ii
yields
$$
  D_A(n_1,\ldots,n_\ell)\geq \Theta n-\theta^* \, .
$$
Thus, using the induction hypothesis,
$$\eqalign{
  |V_{4,n}|
  &{}\leq \sum_{A\in\calQ_+\setminus\widehat\calQ} |r_A|
    \sum_{n_1+\cdots+n_\ell=n-a} |q|^{-\Theta n+\theta^*} c^\ell G^{n-a} \cr
  &{}\leq q^{-\Theta n} n^{\ell^*} \sum_{A\in\calQ_+\setminus\widehat\calQ}
    |r_A| |q|^{\theta^*} c^{\ell^*} G^{n-1} \, .\cr}
$$
Thus, there exists some contant $C_4$ such that
$$
  |V_{4,n}|\leq C_4 c^{\ell^*} G^{n-1} \, .
$$
Combining all these upper bounds, we obtain
$$
  |V_{1,n}+V_{2,n}+V_{3,n}+V_{4,n}|
  \leq (C_1+C_2+C_3+C_4) c^{\ell^*} G^{n-1} \, .
$$
Considering \gnRecA\ for $n\geq p+1$,
$$
  \Bigl|\sum_{A\in\calQ_0} r_A q^{\alpha_1 n}\Bigr| |g_n|
  \leq (C_1+C_2+C_3+C_4) c^{\ell^*} G^{n-1} \, .
$$
Using Lemma \simpleBound,
we then take $G$ to be any number we like greater than $1$ and
$$
  {C_1+C_2+C_3+C_4\over \inf_{n\in\NN}\bigl|\sum_{A\in\calQ_0} 
  r_Aq^{\alpha_1 n}\bigr|}
  c^{\ell^*-1}
$$
to obtain that indeed $|g_n|\leq c\, G^n$.\hfill\qed

\bigskip

It follows from Lemma \gnGeom\ that the generating function 
$g(z)=\sum_{n\geq 0} g_n z^n$ is finite in a neighborhood of the origin. Let us
define
$$
  \tilde g(z)=\sum_{n\geq 0} |g_n| |z|^n \, .
$$
The radius of convergence of $g$ and $\tilde g$ are identical, and
we write $\rho$ for this common radius. We can now prove Theorem \mainTh.

We consider identity \gnRecA, multiply both sides by $z^n$ and sum over $n$.
The left hand side of \gnRecA\ provides the term
$$
  -\sum_{A\in\calQ_0} r_A g(q^{\alpha_1} z)
  = -\sum_{A\in\calQ_0\cap\widehat\calQ} r_A g(z) 
    - \sum_{A\in\calQ_0\setminus\widehat\calQ} r_A g(q^{\alpha_1} z) \, .
$$
Since $\alpha(\calQ_0)=0$ in the assumption of Theorem \mainTh,
if $A$ is in $\calQ_0\setminus\widehat\calQ$, then $\alpha_1\leq -1$, so that
$g(q^{\alpha_1} z)$ has radius of convergence at least $q\rho$.

For $n$ greater than the degree $p$ of $P$, the right hand side 
of \gnRecA\ is decomposed into $V_{1,n}$, $V_{2,n}$, $V_{3,n}$
and $V_{4,n}$ as in the proof of Lemma \gnGeom. Given \gnGeomA, we have
$$
  V_1(z)=\sum_{n\geq 0} V_{1,n}z^n
  = \sum_{A\in\calQ_+\cap\widehat\calQ} r_A s(A) q^{-Ha(a-h)} g_0^{\ell-1} z^a
  g(z) \, .
$$
Thus,
$$
  V_1(z)+\sum_{A\in\calQ_0} r_A g(q^{\alpha_1} z)
  =\calC_{q,g_0}(z)g(z)+\sum_{A\in\calQ_0\setminus\widehat\calQ} r_A g(q^{\alpha_1} z) 
  \, .\eqno{\equa{FinalA}}
$$
Next, the part from $V_{2,n}$ is
$$\eqalign{
  V_2(z) 
  &{}= \sum_{n\geq 0} V_{2,n}z^n\cr
  &{}= \sum_{A\in\calQ_+\cap\widehat\calQ} r_A 
    \sum_{\ss 1\leq i\leq \ell\atop\ss\alpha_i\not=\alpha_\ell} \sum_{n\geq 0}
    q^{-D_A(0,\ldots,n-a,\ldots,0)_i}g_{n-a} z^n g_0^{\ell-1} \, .\cr}
$$
Given \gnGeomB,
$$\eqalign{
  \sum_{n\geq 0} |q|^{-D_A(0,\ldots,n-a,\ldots,0)_i} |g_{n-a}| |z^n|
  &{}\leq C\sum_{n\geq 0} |q|^{-n} |g_{n-a}| |z|^n\cr
  &{}= C |z/q|^a \tilde g(z/q) \, .\cr}
$$
Therefore, the radius of convergence of $V_2(z)$ is at least $q\rho$.

Similarly, the radius of convergence of $V_3(z)=\sum_{n\geq 0} V_{3,n}z^n$
is at least the smaller of the radiuses of convergence of
$$
  V_{A,3}(z)
  =\sum_{n\geq 0} 
  \sum_{\ss n_1+\cdots+n_\ell=n-a\atop\ss \sharp\{i:n_i>0\}\geq 2}
  q^{-D_A(n_1,\ldots,n_\ell)} g_{n_1}\cdots g_{n_\ell} z^n
$$
where $A$ runs over $\calQ_+\cap\widehat\calQ$. Using Lemma \DminCor.i, 
$|V_{A,3}(z)|$ is at most
$$\displaylines{\quad
    \sum_{n\geq 0} 
    \sum_{\ss n_1+\cdots+n_\ell=n-a\atop\ss\sharp\{i:n_i>0\}\geq 2}
    |q|^{-\Theta n+\theta^*} |g_{n_1}|\cdots |g_{n_\ell}| 
    |z|^{n_1+\cdots+n_\ell+a}
    \hfill\cr\noalign{\vskip 3pt}\hfill
  \eqalign{
  {}\leq{}& |q|^{\theta^*} \sum_{n\geq 0} \sum_{n_1+\cdots+n_\ell=n-a} 
    |q^{-\Theta}z|^{n_1}|g_{n_1}|\cdots |q^{-\Theta}z|^{n_\ell}|g_{n_\ell}|
    |q|^{-\Theta a} |z|^a \cr
  {}\leq{}& |q|^{\theta^*} |z|^a\tilde g(z/q^\Theta)^\ell \, .\cr}
  \cr}
$$
Therefore, the radius of convergence of $V_{A,3}$ is at least $q^\Theta\rho$
and so is that of $V_3(z)$. 

Finally, the radius of convergence of $V_4(z)=\sum_{n\geq 0} V_{4,n}z^n$ is at
least the smaller of the radiuses of convergence of
$$
  V_{A,4}(z)=\sum_{n\geq 0} \sum_{n_1+\cdots+n_\ell=n-a} 
  q^{-D_A(n_1,\ldots,n_\ell)} g_{n_1}\cdots g_{n_\ell} z^n
$$
where $A$ is in $\calQ_+\setminus\widehat\calQ$. Using Lemma \DminCor.ii, we obtain
that $|V_{A,4}(z)|$ is at most
$$\displaylines{\quad
  \sum_{n\geq 0} \sum_{n_1+\cdots+n_\ell=n-a} 
    |q|^{-\Theta(n_1+\cdots+n_\ell+a)+\theta^*} |z|^{n_1+\cdots+n_\ell+a}
    |g_{n_1}|\cdots |g_{n_\ell}| 
  \hfill\cr\hfill
  {}\leq |z|^a |q|^{\theta^*-\Theta a} \tilde g(z/q^\Theta)^\ell \, .\qquad\cr}
$$
Therefore, the radius of convergence of $V_4(z)$ is at least $|q|^\Theta\rho$.

Combining \FinalA\ with our estimates on the radius of convergence 
of $V_1(z),\ldots,V_4(z)$, we see that the function
$$
  U(z)=\sum_{A\in\calQ_0\setminus\widehat\calQ} r_A g(q^{\alpha_1} z)
  +V_2(z)+V_3(z)+V_4(z)
  +\sum_{0\leq i\leq p} P_i q^{-d_i} z^i
  \eqno{\equa{FinalC}}
$$
has radius of convergence at least $q^{\min(\Theta,1)}\rho$. Moreover, given
\gnRecA\ and that $g_0=f_0$, we have
$$
  \calC_{q,f_0}(z)g(z)=-U(z) \, .
  \eqno{\equa{FinalB}}
$$
Since $\rho$ is positive and $\calC_{q,f_0}$ is a polynomial in $z$, this 
relation forces $\rho$ to be the smallest modulus of the zeros of the 
crest polynomial. This proves the first assertion of Theorem \mainTh.
Since $U$ has radius of convergence $q^\Theta \rho$, identity
\FinalB\ shows that $\calC_{q,f_0}(z) g(z)$ has removable singularities
in the open disk centered at the origin and of radius $q^\Theta\rho$. 
This proves the second assertion of Theorem \mainTh.

To prove the third assertion, consider the identity \gnRecA. 
If $\sum_{A\in\calQ_0}r_A q^{\alpha_1 n}$ is negative, its left hand side is of
the sign of $g_n$. If the $P_i$ are nonnegative as well as the $r_A$ for \
$A$ in $\calQ_+$, then we can use \gnRecA\ to show inductively that the sign
of $g_n$ is positive whenever that of $g_0$ is.

Finally, if $\calQ_0$ has a unique
element, it must be $(0;0)$ and it belongs to the crest. Thus, \FinalC\
becomes
$$
  U(z)
  = V_2(z)+V_3(z)+V_4(z) +\sum_{0\leq i\leq p} P_i q^{-d_i} z^i \, .
  \eqno{\equa{FinalD}}
$$
Since $(f_n)$ is nonnegative, so is $(g_n)$. From their definition, we
check that $V_{2,n}$, $V_{3,n}$, $V_{4,n}$ and $P_i$, $0\leq i\leq p$,
have all the same sign under all the assumptions of Theorem \mainTh.iii. 
Thus $\bigl([z^n]U(z)\bigr)$ is a sequence of constant sign, which is that 
of $r_A$, $A\in\calQ_+$. The only way for \FinalD\ to vanish is that $(g_n)$
is a sequence whose elements are all $0$, that is $(f_n)=0$ or $f=0$.\hfill\qed

\bigskip

\noindent{\bf Acknowledgements.} We thank Jos\'e Cano and Pedro Fortuny Ayuso
for their help with example 5 in section 4. We also thank Lucia Di Vizio for 
some useful exchanges.


\bigskip

\noindent{\bf References.}
\medskip

{\leftskip=\parindent \parindent=-\parindent
 \par

C.R.\ Adams (1929). On the linear ordinary $q$-difference equations, {\sl Ann.\
Math.}, 30, 195--205.

C.R.\ Adams (1931). Linear $q$-difference equations, {\sl Bull.\ Amer.\ Math.\
Soc.}, 37, 361--400.



J.-P.\ B\'ezivin (1992). Sur les \'equations fonctionnelles aux 
$q$-diff\'erences, {\sl Aequationes Math.}, 43, 159--176.

G.\ Birkhoff (1913). The generalized Riemann problem for linear differential
equations and the allied problems for linear difference and $q$-difference
equations, {\sl Proc.\ Amer.\ Acad.}, 49, 521--568.

J.\ Cano, P.\ Fortuny Ayuso (2012). Power series solutions of nonlinear
$q$-difference equations and the Newton-Puiseux polygon, {\tt arxiv:1209.0295}.

R.D.\ Carmichael (1912). The general theory of linear $q$-difference equations,
{\sl Amer.\ J.\ Math.}, 34, 147--168.

L.\ Di Vizio (2008). An ultrametric version of the Maillet-Malgrange theorem 
for nonlinear $ q$-difference equations, {\sl Proc.\ Amer.\ Math.\ Soc.}, 136, 
2803--2814.

B.\ Drake (2009). Limit of areas under lattice paths, {\sl Discrete Math.},
309, 3936--3953.

Ph.\ Flajolet, R.\ Sedgewick (2009). {\sl Analytic Combinatorics}, Cambridge.

J.\ F\"urlinger, J.\ Hofbauer (1985). $q$-Catalan numbers, {\sl J.\ Comb.\
Th., A}, 248--264.

S.\ Garoufalidis (2004). On the characteristic and deformation varieties of
a knot, {\sl Proceedings of the Casson Fest, Geometry and Topology Monographs},
7, 291--310.

A.M.\ Garsia (1981). A $q$-analogue of the Lagrange inversion formula,
{\sl Houston J.\ Math.}, 7, 205--237.

A.M.\ Garsia, M.\ Haiman (1996). A remarkable $q,r$-Catalan sequence 
and $q$-Lagrange inversion, {\sl J.\ Algebraic Combin.}, 5, 191--244.

I.\ Gessel (1980). A noncommutative generalization and $q$-analog of the 
Lagrange inversion formula, {\sl Trans.\ Amer.\ Math.\ Soc.}, 257, 455-482.

N.\ Joshi (2012). Quicksilver solutions of a $q$-difference first Painlev\'e
equation, {\tt arXiv:1306.5045}.

X.\ Li, C.\ Zhang (2011). Existence of analytic solutions to analytic
nonlinear $q$-difference equations, {\sl J.\ Math.\ Anal.\ Appl.}, 375, 
412--417.

E.\ Maillet (1903). Sur les s\'eries divergentes et les \'equations 
diff\'erentielles, {\sl Ann.\ Sci.\ \'Ec.\ Norm.\ Sup\'er.}, 487--518, 20.

B.\ Malgrange (1989). Sur le th\'eor\`eme de Maillet, {\sl Asymptot.\ Anal.}, 
2, 1--4.

I.G.\ MacDonald (1995). {\sl Symmetric Functions and Hall Polynomials}, Oxford
University Press.

S.\ Nishioka (2010). Transcendence of solutions of $q$-Painlev\'e equation
of type $A_7^{(2)}$, {\sl Aequationes Mathematica}, 79, 1--12.

A.\ Ramani, B.\ Grammaticos (1996). Discrete Painlev\'e equation: coalescence,
limits and degeneratrics, {\sl Phys.\ A}, 228, 160--171.

J.-P.\ Ramis (1992). About the growth of entire functions solutions of linear
algebraic $q$-difference equations, {\sl Ann.\ Fac.\ Sci.\ Toulouse},  6 (1),
53--94.

J.-P.\ Ramis, J.\ Sauloy, C.\ Zhang (2013). {\sl Local Analytic Classification
of $q$-Difference Equations}, preprint.

H.\ Sakai (2001). Rational surfaces associated with affine root systems and
geometry of the Painlev\'e equations, {\sl Comm.\ Math.\ Phys.}, 220, 165--229.

J.\ Sauloy (2000). Syst\`emes aux $q$-diff\'erences singuliers r\'eguliers:
classification, matrice de connexion et monodromie, {\sl Ann.\ Inst.\ Fourier}, 50, 1021--1071.

J.\ Sauloy (2003). Galois theory of Fuchsian $q$-difference equations,
{\sl Ann.\ Sci.\ \'Ecole Norm.\ Sup.}, 36, 925--968.

W.J.\ Trjitzinsky (1938). Theory of non-linear $q$-difference systems,
{\sl Ann.\ Math. Pura Appl.}, 17, 59--106.

C.\ Zhang (1998). Sur un th\'eor\`eme de Maillet-Malgrange pour les \'equations
$q$-diff\'erentielles, {\sl Asymptot.\ Anal.}, 17, 309--314.


C.\ Zhang (2002). A discrete summation for linear $q$-difference equations
with analytic coefficients: general theory and examples, in Braaksma et al.\
ed., {\sl Differential Equations and the Stokes Phenomenon, Proceedings of the 
conference, Groningen, Netherlands, May 28--30, 2001}, World Scientific.

}

\bigskip\bigskip

\setbox1=\vbox{\halign{#\hfil &\hskip 40pt #\hfill\cr
  Ph.\ Barbe                  & W.P.\ McCormick\cr
  90 rue de Vaugirard         & Dept.\ of Statistics \cr
  75006 PARIS                 & University of Georgia \cr
  FRANCE                      & Athens, GA 30602 \cr
  philippe.barbe@math.cnrs.fr & USA \cr
                              & bill@stat.uga.edu \cr}}
\box1

\bigskip\bigskip

\section{Appendix}
We reproduce the {\tt Maple} code which we used in example 5 of section \fixedref{4}.

The first procedure, {\tt getCoeff(R)} takes a polynomial $R$ and solves
for $f_0$ in the equation $R(0,f_0,\ldots,f_0)=0$.

\verbatim@
> getCoeff:=proc(R):
> RootOf(factor(subs(z=0, Y0=c, Y1=c, Y2=c, Y3=c, 
                          Y4=c, Y5=c, R)), c);
> end proc:
@

\bigskip

The next procedure {\tt isReduced(R)} takes a polynomial $R$ and returns 
{\tt true} if the equation given by that polynomial is reduced.

\verbatim@
> isReduced:=proc(R):
> evalb(degree(subs(z=0, R), [Y0,Y1,Y2,Y3,Y4,Y5])
                                              = 1):
> end proc:
@

\bigskip

The next procedure {\tt reduceEquation(R,f0)} takes a polynomial $R$; it 
makes the
change of function $f(z)=f0+zg(z)$ in the equation $R(z,Y_0,\ldots,Y_5)f(z)=0$
and returns the new polynomial that encodes the new equation in $g$.

\verbatim@
> reduceEquation:=proc(R,f0):
> Q := simplify(subs(Y0=f0+z*Y0, Y1=f0+z*r*Y1, 
                     Y2=f0+z*r^2*Y2, Y3=f0+z*r^3*Y3, 
                     Y4=f0+z*r^4*Y4, Y5=f0+z*r^5*Y5,
                     R));
> lz := ldegree(Q,z);
> simplify(Q/z^lz);
> end proc:

@

We input the polynomial in $r$ (not in $q$) accounting for the change
of function $g(z)=f(z^2)$.

\verbatim@
P0 := 4*Y1^4-9*Y0^2*Y1*Y2+2*Y0^3*Y2+z^2*Y0*Y2/r^8
      -z^6*Y0^4*Y5^2-z^6*Y2/r^8-z^6*Y0+z^10

@

Then, we proceeed with the reduction steps.

\verbatim@
> isReduced(P0);
                             false
> g0 := getCoeff(P0);
                               0
> P1 := r^6*reduceEquation(P0, g0):
> isReduced(P1);
                             false
> g1 := getCoeff(P1);

@

At this step the system tells us that $g_1$ it is a root of
$$
  g_1^2\bigl( 1+g_1^2r^8(r-2)(4r-1)+1\bigr) \,.
$$
When $r=2$ or $r=1/4$, the only possibility is $g_1=0$, but in general 
one could choose $g_1$ differently. We follow Cano and Fortuny Ayuso's 
path and take $g_1=0$.

\verbatim@
> g1 := 0:
> P2 := reduceEquation(P1, g1)/r^2:
> isReduced(P2);
                             false
> g2 := getCoeff(P2);
> P3 := reduceEquation(P2, g2)/r^2:
> isReduced(P3);
                             false
> g3 := getCoeff(P3);
                               0
> P4 := reduceEquation(P3, g3)/r^2:
> isReduced(P4);
                             false
> g4 := getCoeff(P4);
                               1
> P5 := reduceEquation(P4, g4)/r^2:
> isReduced(P5);
                             false
> g5 := getCoeff(P5);
                               0
> P6 := reduceEquation(P5, g5)/r^2:
> isReduced(P6);
                             false
> g6 := getCoeff(P6);
                               0
> P7 := reduceEquation(P6, g6)/r^2;
> isReduced(P7);
                             false
> g7 := getCoeff(P7);
                     /   2      6      10     2\
               RootOf\2 r  - 9 r  + 4 r   + _Z /
> alias(a = RootOf(2*r^2-9*r^6+4*r^10+_Z^2)): 
  P8 := reduceEquation(P7, g7);
> isReduced(P8);
                              true
> g8 := getCoeff(P8);
                               0
> P9 := reduceEquation(P8, g8);
> isReduced(P9);
                              true
> g9 := getCoeff(P9);
                               0
> P10 := reduceEquation(P9, g9);
> isReduced(P10);
                              true
> g10 := sort(getCoeff(P10),r);

@

$$
  g_{10}=-{r^2(16r^{11}-9r^{10}-9r^7+2r^6--18r^4+6)\over r^6+1}
$$     

We can check how many monomials are in P10:

\verbatim@

> nops(P10)
                              397

@

Now we have a reduced equation with a solution whose constant term
does not vanish. We calculate $P_0$ and $P_+$.

\verbatim@
> P := P10:
> P0 := subs(z=0, P)-subs(z=0, Y0=0, Y1=0, Y2=0, 
             Y5=0, P):
> PPlus := P-P0-subs(Y0 = 0, Y1 = 0, Y2 = 0, 
                     Y5 = 0, P); nops(PPlus);
                              292
@

This gives $P_0=a(Y_0+r^6Y_2)$, while $P_+$ has $292$ terms.
We then calculate $\alpha(\calQ_0)$.

\verbatim@
> diff(P0, Y5); diff(P0, Y2);
@

We obtain $\partial P_0/ \partial Y_0=0$ and
$\partial P_0/\partial Y_2=r^6a$. Therefore $\alpha(\calQ_0)=2$. We then 
calculate $\alpha(\calQ_+)$.

\verbatim@
> evalb(diff(PPlus, Y5) = 0);
                             false

@

We then have $\alpha(\calQ_+)=5$. We may then check condition
\algoExistenceConditionA. However, this only confirms the result of
Cano and Fortuny Ayuso (2012).

\verbatim@ 
> simplify(subs(z=0,P)/a);

@

We obtain
$$
  Y_0 + 2 r^8 + 16 r^{13} + 6 r^2 - 9 r^9 - 9 r^{12} - 18 r^6 + r^6 Y_2
  \, .
$$
To apply Theorem \mainTh, we need $\alpha(\calQ_0)=0$. This amounts to 
make the change of function $f(z)=g(z/(r^2))$. For the new unkown $g$, 
the equation is obtained by the change of variables 
$(Y_0,Y_1...,Y_5)\rightarrow (Y_{-2},Y_{-1},...,Y_3)$.

\verbatim@
> NP0 := subs(Y0=Ym2, Y1=Ym1, Y2=Y0, Y5=Y3, P0);
> NPPlus := subs(Y0=Ym2, Y1=Ym1, Y2=Y0, Y5=Y3, 
                 PPlus);

@

We then calculate the heights of $R_{-2},\ldots,R_3$. The heights
of $R_{-2}$ and $R_{-1}$ are negative, and cannot be maximal, while that
of $R_0$ vanishes. So we need only the heights of $R_1$ and $R_2$.

\verbatim@
> R1 := subs(Y2=0, Y3=0, NPPlus)-subs(Y1=0, Y2=0, 
  Y3=0, NPPlus): H1:=1/(2*ldegree(R1, z));
                               0
> R2 := subs(Y3=0, NPPlus)-subs(Y2=0, Y3=0, NPPlus):
  H2 := 2/(2*ldegree(R2, z));
                               0
> R3 := NPPlus-subs(Y3=0, NPPlus):
  H3 := 3/(2*ldegree(R3, z));
                               3 
                               --
                               34

@

We obtain $H_1=H_2=0$ and $H_3=3/34$. Thus, the height is  given by $H_3$, 
that is  $3/34$. Since $3/(2H_3)=17$, we need the term in $z^{17}$ in $R_3$.

\verbatim@
> M := coeff(R3, z, 17);

@

This coefficient is $-2r^{64}Y_3$. For this term, $r_A=-2r^{64}$, $s(A)=1$,
$a=17$ and the height is also $17$. To calculate the crest polynomial,
we also need the term involving the $r$-factor $(0;0)$,

\verbatim@

> subs(Ym2=0, Ym1=0, Y1=0, Y2=0, Y3=0, NP0);

@

\noindent
which is $\rho r^6Y_0$. Thus, the crest polynomial is
$$
  \calC_{r,t}(z)=\rho r^6-2r^{64} z^{17} \, .
$$

\bye

\vfill\eject
\bigskip\bigskip

OLD MATERIAL (Bill, don't bother with it; this is for me to remember to rewrite
something seperately and send it to you).

\bigskip\bigskip

One may wonder why we restricted $\calQ$ to be shifting in Theorem \mainTh.
Non-shifting $q$-factors pose existence problems as the following result
shows.\note{We need to exclude the $q$-factor $(0;0)$ from $\calQ$ everywhere
since this yields another $f(z)$ which could cancel with the left hand side}

\Proposition{\label{noSolution}
  Assume that $P$ is a polynomial with nonnegative coefficients, that 
  each $r_A$ in \qEq\ is nonnegative, and that $\calQ$ contains a nonshifting
  $q$-factor with $\alpha_\ell$ positive. If $q>1$, then \qEq\ has no formal 
  power series solution except
  perhaps a polynomial.
}

\bigskip

\Proof Assume that $f$ solves \qEq\ and is not a polynomial. 
Let $A=(0;\alpha_1,\ldots,\alpha_\ell)$\note{item 37}
be a nonshifting $q$-factor in $\calQ$ with $\alpha_\ell$ positive.
We will see in the proof of Theorem
\mainTh, that if the coefficients of $P$ are positive as well as the $r_A$, 
then a formal solution has nonnegative coefficients. Applying $[z^n]$ to both
sides of \qEq, we then obtain
$$\eqalign{
  [z^n]f(z)
  &{}\geq r_A [z^n]Af \cr
  &{}=r_A\sum_{n_1+\cdots+n_\ell=n} q^{\alpha_1n_1+\cdots+\alpha_\ell n_\ell}
    f_{n_1}\cdots f_{n_\ell} \, .\cr
  }
$$
Consider only the tuple $(n_1,\ldots,n_\ell)=(0,\ldots,0,n)$ to obtain
$$
  f_n\geq r_A q^{\alpha_\ell n} f_n \ ,.
$$
Since $\alpha_\ell$ does not vanish and $q>1$, this relation cannot hold
for $n$ large enough.\hfill\qed

\bigskip

\bigskip

\subsection{Reduction}

The following theorem is quite obvious if one examine its meaning on a 
few concrete examples.

\Theorem{\label{reductionTh}
  If $f$ solves \qEq, then exists a formal power series $g$ such that
  $f=f_0+zg(q^{-\alpha(\calQ_0)}z)$ and $g$ solves a reduced $q$-algebraic
  equation.
}

\bigskip

To prove this theorem, we introduce the following definition which will
make some bookeeping much easier. We 
write $(\alpha_1,\ldots,\alpha_\ell)_{[k]}$ the set of
subsequence of $\{\, 1,2,\ldots,k\,\}$  of length $k$. There 
are ${\ell\choose k}$ such subsequences.
This set is empty if $k>\ell$. If $k$ or $\ell$ is $0$, we agree that this
set is empty.

\Definition{\label{qDiff}
  Let $A=(0;\alpha_1,\ldots,\alpha_\ell)$ be a nonshifting $q$-factor 
  and let $t$ be a complex number. The symbolic $k$-th differential of $A$ 
  in  the direction $t$ is the linear combination of $q$-factors $D_t^kA$ 
  defined by
  $$
    D_t^k
    =k!t^{\ell-k} \sum_{\theta\in (\alpha_1,\ldots,\alpha_\ell)_[k]} 
     q^{\theta_1+\cdots+\theta_k} f(q^{\theta_1}z)\cdots f(q^{\theta_k}z) \, .
  $$
}

Note that $D_t^1A$ is the linear operator $f\mapsto t^\ell\sum_{1\leq i\leq\ell}
q^{\alpha_i}f(q^{\alpha_i}z)$. The following lemma is the motivation to introduce
definition \qDiff.

\Lemma{\label{qTaylor}
  For any nonshifting $q$-factor $A$ and any complex number $t$,
  $$
    A\bigl( t+zg(z)\bigr) =\sum_{0\leq k\leq \ell} z^k {D_t^kA\over k!}g(z)\,.
  $$
}

\Proof We have
$$
  A\bigl( t+zh(z)\bigr)
  =\prod_{1\leq i\leq\ell} \bigl( t+q^{\alpha_i}z h(q^{\alpha_i}z)\bigr) \, .
$$
We then expand this product on the power of $z$.\hfill\qed

\bigskip

One more definition will be handy.

\Definition{\label{nonshiftingPart}
  Let $A=(a;\alpha_1,\ldots,\alpha_\ell)$ be a $q$-factor. Its nonshifting
  part is the $q$-factor $(0;\alpha_1,\ldots,\alpha_\ell)$. 
}

\bigskip

\noindent{\bf Proof of Theorem \reductionTh.} Apply $[z^0]$ to both sides
of \qEq\ to obtain
$$
  \sum_{A\in\calQ_0}r_A f_0^\ell=0 \, .
  \eqno{\equa{reductionA}}
$$
Let $h(z)=g(q^{-\alpha(\calQ_0)}z)$. Using Lemma \qTaylor,
$$\eqalign{
  \sum_{A\in\calQ} r_QA\bigl(f_0+zh(z)\bigr)
  &{}= \sum_{A\in\calA}r_A z^a A_0\bigl( f_0+zh(z)\bigr) \cr
  &{}=\sum_{A\in\calQ} r_A z^a \Bigl( f_0^\ell +\sum_{k\geq 1} z^k 
    {D_t^kA_0\over k!} h(z)\Bigr) \, . \cr}
$$
Given \reductionA, this is the sum of two terms,
$$
  \sum_{A\in\calA\setminus\calQ_0}r_A z^a f_o^\ell \,,
$$
and
$$
  \sum_{A\in\calQ} r_A z^a \sum_{k\geq 1} z^k {D_t^kA_0\over k!} h(z) \, .
$$
Both terms are in the ideal generated by $z$. Therefore, we can rewrite
\qEq\ as
$$
  \sum_{A\in\calQ\setminus\calQ_0} r_A z^{a-1} f_0^\ell
  +\sum_{A\in\calQ} r_A z^a\sum_{k\geq 0} z^k {D_t^{k+1}\over k!} A_0 h(z)
  = 0 \, .
$$
In this new $q$-functional equation, the nonshifting $q$-factors are
$$
  \sum_{\ss A\in\calQ\setminus\calQ_0\atop\ss a=1} r_A z^{a-1} f_0^\ell \, ,
$$
which is a  polynomial in $z$ and
$$
  \sum_{A\in\calQ_0} r_A {D_{f_0}^1 A_0\over k!} h(z) \, .
$$
Since $D_{f_0}^1$ is linear, the nonshifting $q$-factors are linear.\note{except for the polynomial terms}

Recall that $h(z)=\sigma^{-\alpha(\calQ_0)}g(z)$. 
If $A=(a;\alpha_1,\ldots,\alpha_\ell)$ is a $q$-factor, then 
$$
  A\sigma^{-p}=(a;\alpha_1-p,\ldots,\alpha_\ell-p)
$$
and the result follows.\hfill\qed

\bigskip

~\ \note{keep what follows?}

\Lemma{\label{normalization}
  Suppose that $f$ satisfies \qEq.

  \medskip
  (i) Assume that $P(z)=z^k Q(z)$ for some polynomial $Q$ and some 
  positive $k$. For a $q$-factor $A=(a;\alpha_1,\ldots,\alpha_\ell$, 
  set $A_{[k]}=(a+(\ell-1)k; \alpha_1,\ldots,\alpha_\ell)$. We then have 
  $f(z)=z^kg(z)$ where $g$ solves
  $$
    \sum_{A\in\calQ_0} r_A q^{k(\alpha_1+\cdots+\alpha_\ell)} A_{[k]}g(z)
    =Q(z)+\sum_{A\in\calQ} r_A q^{k(\alpha_1+\cdots+\alpha_\ell)} 
    A_{[k]} g(z) \, .\eqno{\equa{normalizationA}}
  $$

  \noindent (ii) For any complex number $\lambda\not=0$ the formal power
  series  $g(z)=\lambda f(z)$ where $g(z)$ solves
  $$
    \sum_{A\in\calQ_0} g(z)
    =\lambda^{-1}P(z) +\sum_{A\in\calQ}r_A \lambda^{\ell-1} Ag(z) \, . 
  $$
}

Note that if $A$ is nonshifting, then $A_{[k]}$ is shifting if and only 
if $\ell>1$. Therefore, the nonshifting $q$-factors in the left hand side 
of \normalizationA\ come only from those of the form $(0;\alpha_1)$ 
in $\calQ_0$, all the others becoming shifting under the transformation
$A\mapsto A_{[k]}$.

\bigskip

\Proof Consider a $q$-factor $(a;\alpha_1,\ldots,\alpha_\ell)$. 

\noindent{\it (i).} We have
$$\displaylines{\qquad
  z^{-k} (a;\alpha_1,\ldots,\alpha_\ell)\bigl(z^kg(z)\bigr)
  \hfill\cr\noalign{\vskip 5pt}\hfill
  \eqalign{
  {}={}&z^{-k} z^a\prod_{1\leq i\leq \ell} (q^{\alpha_i}z)^k g(q^{\alpha_i}z) \cr
  {}={}&q^{(\alpha_1+\cdots+\alpha_\ell)k} 
        \bigl(a+k(\ell-1);\alpha_1,\ldots,\alpha_\ell\big) g(z) \, . \cr}
  \qquad\cr}
$$
The first assertion of the lemma follows.

\noindent{\it (ii).} For any real number $\lambda\not=0$, 
$$
  \lambda^{-1}(a;\alpha_1,\ldots,\alpha_\ell)\bigl(\lambda g(z)\bigr)
  = \lambda^{\ell-1} (a;\alpha_1,\ldots,\alpha_\ell) g(z) \, .
$$
The second assertion follows.\hfill\qed

\bigskip

\vfill\eject

\bye